\documentclass{amsart}
\usepackage[utf8]{inputenc}
\usepackage{amsmath}
\usepackage{amssymb}
\usepackage{mathrsfs}
\usepackage{graphicx}

\usepackage[colorlinks]{hyperref}
\usepackage[svgnames, dvipsnames, usenames]{xcolor}
\hypersetup{urlcolor = RedViolet, linkcolor = RoyalBlue, citecolor = ForestGreen}

\usepackage{tikz}

\newtheorem{thmx}{Theorem}

\newtheorem{theorem}{Theorem}[section]
\newtheorem{corollary}[theorem]{Corollary}
\newtheorem{lemma}[theorem]{Lemma}
\newtheorem{proposition}[theorem]{Proposition}
\theoremstyle{definition}
\newtheorem{definition}{Definition}[section]
\newtheorem{remark}{Remark}[section]
\newtheorem{question}{Question}[section]

\newcommand{\A}{\mathcal{A}}
\newcommand{\C}{\mathbb{C}}
\newcommand{\Cat}{\mathscr{C}}
\newcommand{\F}{\mathcal{F}}
\newcommand{\N}{\mathbb{N}}

\newcommand{\T}{\mathbb{T}}
\newcommand{\Z}{\mathbb{Z}}

\newcommand{\Hopf}{\mathsf{Hopf}}
\newcommand{\AffCat}{\mathsf{Aff}}

\newcommand{\tilstar}{\mathbin{\tilde *}}
\newcommand{\Ctrans}{^{\mathsf{T}}}
\newcommand{\sync}{^{\mathrm{s}}}
\DeclareMathOperator{\Aut}{Aut}
\DeclareMathOperator{\Aff}{Aff}
\DeclareMathOperator{\Tr}{Tr}

\usepackage{tikz}
\usepackage{strdiag}
\usepackage{fancyhdr}
\title{Complex Hadamard Matrices -- Quantum Symmetries, Equivalence and Non-Local Games}
\author[Brannan]{Michael Brannan}
\address{Department of Pure Mathematics and Institute for Quantum Computing, University of Waterloo. 200 University Ave W, Waterloo, ON N2L 3G1, Canada}
\email{michael.brannan@uwaterloo.ca}

\author[Gromada]{Daniel Gromada}
\address{Czech Technical University in Prague, Faculty of Electrical Engineering, Department of Mathematics, Technická 2, 166 27 Praha 6, Czechia}

\email{gromadan@fel.cvut.cz}

\author[Hernández Palomares]{Roberto Hernández Palomares}
\address{Department of Mathematics, College of Arts and Sciences, The Ohio State University. 100 Math Tower 231 West 18th Avenue Columbus, OH 43210-1174 USA}
\email{hernandezpalomares.1@osu.edu}

\author[Priebe]{Nicholas Priebe}
\email{nicholas.priebe@uwaterloo.ca}
\date{\today}

\begin{document}

\maketitle

\begin{abstract}
We consider quantum group generalizations of the action of monomial matrices on complex Hadamard matrices.  This gives rise to various notions of quantum symmetries and quantum equivalences of Hadamard matrices. We show that if one acts by a certain largest monomial quantum group, then all Hadamard matrices of a given size become quantum equivalent.  Taking a more restrictive quantization leads to a notion of  $s$-quantum equivalence. We exhibit examples of Butson matrices of the same size and order that are not $s$-quantum equivalent for any choice of $s$.  We also show that $s$-quantum equivalence is operationally modeled by a synchronous  non-local ``Hadamard equivalence'' game.  Our methods are largely graphical calculus based, using categories generated by complementary spiders.  We use these same tools to study quantum affine equivalence of quantum groups, and prove that all finite quantum groups of the same size are quantum affinely equivalent.  We also provide a short graphical proof of a result of Kasprzak--So\l tan--Woronowicz asserting that  quantum symmetries of finite quantum groups must be classical.             
\end{abstract}

\section{Introduction}

A (complex) Hadamard matrix is a square matrix $H \in M_n(\C)$ whose entries have modulus one and whose rows (equivalently, columns) are pairwise orthogonal. %Thus, after renormalization, a complex Hadamard matrix is precisely the unitary transition matrix between two mutually unbiased orthonormal bases of $\C^n$.
Hadamard matrices  are ubiquitous within  mathematics -- appearing in  various guises in  combinatorics, operator algebras, representation theory, computer science, and (quantum) information theory. See  \cite{Hor07} for applications in computer science and information theory,   \cite{CM22, GHad} for connections to quantum groups and  quantum information, and  \cite{Po83, Nic10} for applications to operator algebras. 

 The prototypical example of a Hadamard matrix is the Fourier transform matrix $\mathcal F_n = [e^{\frac{2\pi i kl}{n}}]_{k,l=0}^{n-1}$ associated to  the finite abelian group $\mathbb Z_n$.  Fourier matrices are special instances of Butson Hadamard matrices -- Hadamard matrices whose entries are  restricted to be roots of unity.  Within combinatorics, a particularly special role is played by the {\it real} Hadamard matrices -- i.e., those whose entries belong to $\{\pm 1\}$.  It is well-known that a real Hadamard matrix  $H \in M_n(\C)$ can exist only when $n=2$ or $n$ is divisible by $4$.  The famous {\it Hadamard conjecture} asserts that this necessary condition is also sufficient.  At present, the Hadamard conjecture is open and has only been verified for all $n <  668$ (the case $n=428$ being unknown until 2005 \cite{KhTa05}).\footnote{On August 12, 2026 Levent Alp\"{o}ge, Saul Reynolds-Haertle and Philippe Voinov announced an AI (Claude LLM)-assisted verification of the Hadamard conjecture up to $n=2000$.}  Of course, the existence of Fourier matrices shows that complex Hadamard matrices exist in all dimensions.   

Given a complex Hadamard matrix $H \in M_n(\C)$, one can transform it into another Hadamard matrix $H'\in M_n(\C)$ by  permuting rows and columns and multiplying rows and columns by complex phases.  This gives rise to a natural equivalence relation on Hadamard matrices of a given size, and a fundamental problem in the theory of Hadamard matrices is to classify Hadamard matrices modulo this equivalence relation.  The above natural operations on Hadamard matrices are group-theoretically realized by left and right multiplications by the group of   monomial matrices. For Butson  matrices whose entries are $s$-th roots of unity, one can assume the monomial matrices have entries of the same type.  Thus the relevant symmetries are implemented by left and right multiplications by elements from the  the (finite) complex reflection group $H_n^s$ of complex permutation matrices with values in the group $\hat{\mathbb Z}_s$ of $s$-th roots of unity.  In other words, symmerties and equivalences of size $n$ and  order $s$ Butson matrices are governed by actions of the wreath product groups
\[
H_n^s = \hat{\mathbb{Z}}_s \wr S_n := (\hat{\mathbb{Z}}_s)^n \rtimes S_n.
\]
The purpose of this article is to study quantizations of these group actions on Hadamard matrices -- replacing $H_n^s$ with a genuine quantum group, and studying the corresponding (potentially larger) {\it quantum} equivalence classes of Hadamard matrices.  

Our starting point is the theory of quantum symmetries of real Hadamard matrices developed in earlier work of the second author \cite{GHad}. In the real case, the classical monomial symmetry group is the hyperoctahedral group $H_n:=H_n^s$ at $s=2$, and its natural quantum analogue is the hyperoctahedral quantum group $H_n^+$ \cite{BBC}. This leads one to associate to any real Hadamard its {\it quantum automorphism group} $\Aut^+(H) := H_n^+ \cap H H_n^+ H^{-1}$.  More generally, to any pair $H,H'$ of real Hadamard matrices of the same size, one associates  a compact quantum space $H_n^+ \cap H H_n^+ {H'}^{-1}$ encoding the (a priori possibly empty) space of quantum equivalences between $H$ and $H'$.  This space of quantum equivalences is modeled by a certain linking $*$-algebra $\mathcal A(H,H')$, which is non-zero iff $H$ and $H'$ are quantum equivalent iff the  compact quantum groups $\Aut^+(H)$ and $\Aut^+(H')$ are canonically monoidally equivalent.  A   striking result in \cite{GHad} is that any two real Hadamard matrices of the same size are quantum equivalent, even when they are not classically equivalent.  Interpreting this result in terms of the bipartite Hadamard graphs $\Gamma_H, \Gamma_{H'}$ associated to $H,H'$ \cite{McK79} (see also Section \ref{sec:graphs}), one obtains that any pair of Hadamard graphs associated to real Hadamard matrices are quantum isomorphic in the quantum information-theoretic sense \cite{AMR19}.  This quantum isomorphism result for graphs was independently obtained by Chan--Martin \cite{CM22} using different methods.   

The above results for real Hadamard matrices lead to many interesting questions about what happens in the complex case. For example, could it be that that all Butson matrices of a given size $n$ and order $s$ are quantum equivalent in a way that generalizes the real case? We show that the  passage from real to complex Hadamard matrices introduces several new phenomena. There is no longer a unique natural quantum analogue of the groups of monomial matrices. For matrices of fixed size and order $s$, a natural choice is to work with the complex quantum reflection groups $H_n^{s+} = \hat{\mathbb Z}_s \wr_* S_n^+$ \cite{BV09}.  Quantizing symmetries with arbitrary complex phases leads to two distinct  quantum groups: The free wreath product  $H_n^{\infty+} = \mathbb T \wr_* S_n^+$ and the free complexification of $H_n^+$, denoted $K_n^+ := H_n^+ \widetilde{*}\, \mathbb{T}.$
These distinct choices of quantum groups and their actions on Hadamard matrices can be used to  produce three corresponding notions of quantum equivalence of Hadamard matrices, which we call $s$-quantum equivalence, strong quantum equivalence, and weak quantum equivalence. By definition these quantum equivalences (QE) satisfy the following implications 
\[
s\text{-QE}
\Longrightarrow
\text{strong QE}
\Longrightarrow
\text{weak QE}.
\]
The above notions of quantum equivalence are defined analogously to the real case: For complex Hadamard matrices $H, H'$ of size $n$, we define the quantum automorphism group of $H$ relative to one of the above quantum reflection type groups $G = H_n^{s+}, H_n^{\infty +}, K_n^+$ to be
\[
\Aut_G^+(H)
=
G \cap H G H^{-1},
\]
and we also define a universal linking $*$-algebra $\mathcal A_G(H,H')$. The non-vanishing of this algebra is, by definition, equivalent to quantum equivalence of $H$ and $H'$ relative to the quantum group $G$. These algebras are bi-Galois objects connecting the corresponding quantum automorphism groups. Consequently, $H$ and $H'$ are $G$-quantum equivalent precisely when the representation categories of $\Aut_G^+(H)$ and $\Aut_G^+(H')$ are related by a canonical unitary monoidal equivalence preserving the distinguished generating morphisms. We describe these representation categories using the graphical calculus of complementary spiders.          The graphical calculus description of these representation categories allows for diagrammatic methods in place of direct calculations in universal $*$-algebras. Our first main  result  is a positive one:

\begin{thmx}[See Theorem \ref{T.weak}]
All complex Hadamard matrices of a fixed size are weakly quantum equivalent.
\end{thmx}

  Thus, at the level of the largest quantum symmetry group $K_n^+$ considered here, the phenomenon previously observed for real Hadamard matrices persists. As in the real case, this result is obtained by showing that the abstract diagrammatic category generated by relevant complementary families of spiders is pure, in the sense that every closed diagram reduces to a scalar depending only on $n$. It follows that all fibre functors associated with complex Hadamard matrices of the same size have the same kernel.

Strong quantum equivalence and $s$-quantum equivalence turn out to behave quite differently. We exhibit  families of Hadamard matrices of the same size that are not strongly quantum equivalent. We show in Proposition \ref{P.Fourier} that suitable (tensor product) Fourier matrices of the same size are not strongly quantum equivalent. We also exhibit some  examples among Butson matrices of size $14$ and  prime order $7$ in Proposition \ref{P.H7}.  

A second theme of the paper is the interpretation of quantum Hadamard equivalence through the lens of non-local games. For Butson  matrices $H,H'$ of size $n$ and order $s$, we define a synchronous Hadamard equivalence game, in which two non-communicating players attempt to demonstrate an $s$-monomial equivalence between $H$ and $H'$.  We prove that the game algebra in the sense of \cite{HMPS19} is naturally $*$-isomorphic to the linking algebra
$\mathcal A_{H_n^{s+}}(H,H')$
associated to $s$-quantum equivalence.  
We establish these results through a more general theory of monomial isomorphism games for arbitrary matrices introduced in Section \ref{sec:games}.  For real Hadamard matrices, results in \cite{GHad} show that that our Hadamard equivalence game can be reformulated in terms of a graph isomorphism game between two bipartite graphs $\Gamma_H, \Gamma_{H'}$ associated to our Hadamard matrices
\cite{McK79}.  For higher order Butson matrices, we explain how these graphs can still be constructed, but on the other hand show that in general the (quantum) symmetries of these graphs $\Gamma_H$ can be strictly larger than those of $H$.   Thus, the Hadamard equivalence game is, in general, distinct from a graph isomorphism game when the matrices involved are complex.  Nevertheless, we show that the quantum space of  quantum $s$-equivalences of Butson Hadamard matrices embeds into the quantum space of quantum isomorphisms between the associated graphs (cf. Theorem \ref{th:embedding} and Corollary \ref{cor:strategy-transfer}). This gives a complex analogue of one direction of the  bijective correspondence one has for real Hadamard matrices. 

The final part of the paper considers categories generated by pairs of complementary spiders via fibre functors associated to finite quantum groups. Working with finite quantum groups from the perspective of this graphical calculus, we are able to provide very short graphical proof of the following result originally proved in \cite{KSW15}:

\begin{thmx}[\cite{KSW15} and Corollary \ref{cor:noqs}]
Finite quantum groups only admit classical symmetries. 
\end{thmx}
That is, for any finite quantum group $G$, its quantum automorphism group (preserving both the algebra and coalgebra structures) is automatically a classical group.

We also introduce a quantum version of affine equivalence for finite quantum groups. Classically, two finite groups are affinely equivalent if and only if are isomorphic. In contrast, we prove the following striking quantum variant of this result.

\begin{thmx}
 All finite quantum groups of a fixed dimension are mutually quantum affine equivalent.    
\end{thmx}

The above result is reminiscent of the collapse of equivalence classes when one passes from classical to quantum equivalences of real Hadamard matrices.  We leave it for future work to determine a physical (or operational) interpretation of these quantum affine equivalences in terms of, say,  entanglement.  

The remainder of the paper is organized as follows. In Section \ref{sec:prelims} we review Hadamard matrices, their monomial equivalence relations, and their classical automorphism groups. Section \ref{sec:diagrams} introduces the relevant quantum reflection groups, the associated notions of quantum equivalence, and their categorical characterization through linking algebras and monoidal equivalences. Section \ref{sec:examples} proves the weak quantum equivalence theorem and presents examples separating the stronger equivalence relations. Sections \ref{sec:games} and \ref{sec:graphs}
deal with non-local games and Hadamard graphs respectively, and finally Section \ref{sec:qsymQG} deals with (affine)  symmetries of finite quantum groups. 

\subsection*{Acknowledgements} The authors thank Ada Chan, Alex Chirvasitu and Chris Godsil for insightful conversations.    This mathematical research was conducted by humans and was supported by an NSERC Discovery Grant.

\section{Preliminaries} \label{sec:prelims}

\subsection{Notation}

Throughout the article, $A^\dag$ denotes the conjugate transpose of a matrix $A \in M_n(\C)$ (and more generally, the adjoint of a bounded linear map between Hilbert spaces), $A^*$  denotes the conjugate matrix, $\bullet$ denotes the entrywise Schur product of matrices, and $J$ denotes the all-ones-matrix (the unit of $\bullet$).  Thus if $A = [A_{ij}], B = [B_{ij}] \in M_n(\mathbb C)$, then 
\[(A^\dag)_{ij} =\overline{A_{ji}}, \quad(A^*)_{ij} = \overline{A_{ij}}, \qquad (A \bullet B)_{ij} = A_{ij}B_{ij}, \quad J_{ij} = 1.\]
At times we will also apply the Schur product $\bullet$ to operator valued matrices $A,B$ in the obvious way. For $s \in \N$, we denote by $\Z_s$ the cyclic group of order $s$, and its Pontryagin dual group by $\hat \Z_s$.  We regard $\hat\Z_s  = \{\omega^k\}_{k=0}^{s-1}$, where $\omega \in \mathbb T$ is a primitive $s$th root of unity.
\subsection{Hadamard matrices}

In this work, we will deal with \emph{complex} Hadamard matrices unless stated otherwise. That is, we take the following definition:

\begin{definition}
\label{D.Hadamard}
    A \emph{Hadamard matrix} of size $n$ is any matrix $H \in M_n(\C)$ that
    \begin{enumerate}
        \item has complex unit entries, i.e. $H\bullet H^* =J,$
        \item is up to normalization unitary, i.e. $HH^\dag=nI = H^\dag H.$
    \end{enumerate}
    We say that a complex Hadamard matrix $H$ is of order $s\in\N$ if its entries are $s$-th roots of unity, i.e. $H^{\bullet s}=J$ (often called Butson type as they first appeared in \cite{But62}). In particular, $H$ is a \emph{real} Hadamard matrix if it is of order two.
\end{definition}

Note that in addition to arising in varius parts of pure mathematics, Hadamard matrices are also interesting for quantum computation as they stand for transformations between mutually unbiased bases: Let $V$ be a Hilbert space and $(e_i)$, $(f_j)$ a pair of orthonormal  bases. We call them \emph{mutually unbiased} if $|\langle e_i,f_j\rangle|=\frac{1}{\sqrt{n}}$ for every $i,j$. Denoting by $U$ the unitary transformation matrix defined by $Ue_i = f_i$, the condition of $(e_i), (f_i)$ being mutually unbiased exactly translates to $H = \sqrt{n}U$ being a Hadamard matrix.   

\subsection{Hadamard matrix equivalence}

There are several notions of equivalence for Hadamard matrices all generalizing the only common one for the real case.  In this paper, a \emph{monomial matrix} is a matrix $P \in M_n(\mathbb C)$ whose entries either have modulus 1 or 0, and with exactly one non-zero entry in each row and each column.  Note that a monomial matrix with entries in $\{0,1\}$ is an  $n\times n$ permutation matrix.  

\begin{definition}
    Two Hadamard matrices $H_1$ and $H_2$ of size $n$ are called \emph{monomially equivalent} if there is a pair of monomial matrices $P,Q \in M_n(\T)$ such that $H_1=PH_2Q^{-1}$. We say that they are {\it $s$-monomially equivalent} if there exist such monomial matrices whose non-zero entries are $s$-th roots of unity.

    Equivalently, $H_1$ and $H_2$ are monomially equivalent iff $H_1$ can be obtained from $H_2$ by performing a sequence of the following allowed  operations: 
    \begin{enumerate}
        \item
        Permuting  rows and/or columns of a matrix. 
        \item Multiplying the rows and/or columns  of a matrix by complex numbers of modulus $1$.  
    \end{enumerate}
    In other words, $H_1, H_2$ are monomially equivalent iff there exist permutations $\pi, \sigma \in S_n$ (the symmetric group on $n$ letters), and $z_i, w_i \in \mathbb T$ ($1 \le i \le n$) such that $H_{1,ij} = z_iw_j^{-1} H_{2,\pi(i)\sigma(j)}$ for all $1 
    \le i,j \le n$.
\end{definition}

\begin{remark}
    It is worth pointing out that if $H_1, H_2$ are both order $s$ Butson matrices, then they are  monomially equivalent iff they are $s$-monomially equivalent.  To see this claim, suppose  $H_1, H_2$ are monomially equivalent.  Let $\pi, \sigma \in S_n$ and $w_i, z_i \in \T$ be as in the previous definition. Since $H_1,H_2$ are both order $s$, it follows that $z_iw_j^{-1}  \in \hat \Z_s$ for all $i,j$.  Define $\tilde z_i = z_iz_1^{-1} = (z_i w_1^{-1})(w_1 z_1^{-1} )\in \hat \Z_s$ and $\tilde w_i = w_iz_1^{-1} \in \hat \Z_s$.  Then for $1 \le i,j \le n$, we have \[H_{1,ij} = z_iw_j^{-1} H_{2,\pi(i)\sigma(j)} = (\tilde{z_i}z_1)(z_1 \tilde{w_j})^{-1}H_{2,\pi(i)\sigma(j)}  = \tilde z_i\tilde w_j^{-1} H_{2,\pi(i)\sigma(j)}, \] so $H_1, H_2$ are $s$-monomially equivalent.
\end{remark}
If, in addition to the transformations allowed by monomial equivalence, we allow entrywise automorphisms of $\hat \Z_s$, we obtain the notion of Hadamard equivalence.  
\begin{definition}
    Two Hadamard matrices $H_1$ and $H_2$ of order $s$ are called \emph{Hadamard equivalent} if there is a pair of monomial matrices $P$, $Q$ with entries in $s$-th roots of unity and an automorphism $\phi$ of $\hat \Z_s$ such that $H_1=P\phi(H_2)Q^{-1}$. Here, the automorphism should be applied entrywise.
\end{definition}  

We denote by $H_n^s$ the group of monomial matrices of order $s \le \infty$. $H_n^s$ is often called the complex reflection group size $n$ and order $s$.  $H_n^s$ can be identified with the wreath product $H_n^s \cong \hat \Z_s \wr S_n:= (\hat \Z_s)^n \rtimes S_n$, where $S_n$ acts on $(\hat \Z_s)^n$ by permuting factors. When $s=\infty$, we take $\hat \Z_s = \T$.  

Denote by $\text{Had}_s(n)\subseteq M_n(\C)$  the set of order $s$ Hadamard matrices of size $n$. There is a natural action 
\[H_n^s \times H_n^s \curvearrowright \mathrm{Had}_s(n); \qquad H \mapsto PHQ^{-1} \quad (H \in \mathrm{Had}_s(n), \ P,Q \in H_n^s).\]
The monomial equivalence class of $H$ is precisely the orbit of $H$ under this group action. 

Using this setup, one can also define the {\it (order $s$ monomial) automorphism group of $H \in \mathrm{Had}_s(n)$},  $\Aut_{H_n^s}(H)$, to be the stabilizer subgroup of $H$ under the above action.  Equivalently,
\begin{align*}\Aut_{H_n^s}(H) &= \{P, Q \in H_n^s: H = PHQ^{-1} \}\\
&= \{P, Q \in H_n^s: H^{-1}PH = Q\} \\
&=H_n^s \cap( H^{-1}H_n^s H).\end{align*}

One can also describe the notion of Hadamard equivalence group theoretically.  One just has to replace the complex reflection groups $H_n^s = (\hat \Z_s)^n \rtimes S_n$ with the semidirect products $(\hat \Z_s)^n \rtimes (S_n \times \text{Aut}(\hat \Z_s))$, where $S_n \times \text{Aut}(\hat \Z_s)$ acts on $(\hat \Z_s)^n$ via 
\[(\sigma, \phi) \cdot (z_i)_{i=1}^n  = (\phi(z_{\sigma(i)}))_{i=1}^n, \qquad \sigma \in S_n, \phi \in \Aut(\hat \Z_s).\]

In the sequel, we will mainly be interested in quantising the notion of monomial equivalence for Hadamard matrices, so our focus will be on the complex reflection groups $H_n^s$ and their quantum group analogues.  

\section{Quantum symmetries of complex Hadamard matrices} \label{sec:diagrams}

The quantum symmetries that arise in this work are C$^*$-algebraic compact quantum groups \cite{Wor95}.  We describe these objects both as unital C$^*$-algebras (or Hopf $*$-algebras) viewed as quantized algebras of continuous functions on a compact group, and also in terms of their (concrete) representation categories, via Tannaka--Krein reconstruction methods \cite{Wor88, NT13, TW17}.     

\subsection{Categories and diagrams}

Let $V \simeq\C^n$ be a Hilbert space of dimension $n$, and let $V^* \simeq \C^n$ denote the dual space. Denote by $(e_i)$ a fixed orthonormal basis for $V$, and let  $(e_i^*)$ be the associated dual basis of $V^*$. We consider the canonical pairings $R\colon V^*\otimes V\to\C$ and $R^*\colon V\otimes V^*\to\C$ given by $R(e_i^*\otimes e_j)=\delta_{ij}=R^*(e_i\otimes e_j^*)$. We extend the operation $e_i\mapsto e_i^*$ antilinearly to obtain an antilinear isomorphism $*\colon V\to V^*$.

In this work, a \emph{concrete representation category} is a rigid monoidal $\dag$-category with objects being all possible tensor products of $V$ and $V^*$ containing $R$ and $R^*$ (making $V$ and $V^*$ duals of each other in the monoidal category sense) \cite{NT13, EGN}. That is, it is a family of linear maps between such tensor products of $V$ and $V^*$ such that it is closed under linear combinations (whenever possible), compositions (whenever possible), tensor products (denoted by $\otimes$), and taking adjoints (here denoted by $\dag$). If $A,B,C,\dots$ are such linear maps, we denote by $\langle A,B,C, \ldots \rangle$ the smallest monoidal $\dag$-category containing them together with $R$, $R^*$.

We will use the following diagrammatic approach to define concrete representation categories. We denote (the identity on) $V$ by $\Fid$, (the identity on) $V^*$ by $\Bid$. We denote the dual morphisms by cups and caps:
$$R=\Rcup,\quad R^*=\Lcup,\quad R^\dag=\Rcap,\quad R\Ctrans:=(R^*)^\dag=\Lcap.$$

In addition, we define the \emph{spiders}
$$\Diagram{\Dmor{bcirc}[<.</>.>] (0,.5)}=\sum_i e_i\otimes\cdots\otimes e_i\otimes e_i^\dag\otimes\cdots\otimes e_i^\dag,$$
where we can also flip any arrows replacing $e_i$ by $e_i^*$ on the right-hand-side. In particular, $\Diagram{\Dmor{bcirc}[</<] (0,0.5)}$ represents the canonical isomorphism $V\to V^*$ mapping $e_i\mapsto e_i^*$.  If this particular map is contained in our category, distinguishing the direction of the arrows is pointless as this morphism can be used to flip the arrows arbitrarily. (Although the arrows may still actually be useful for some computations.)  In our categories, composition of (composable) maps is graphically represented by vertical concatenation, and diagrams are read from bottom (inputs) to top (outputs).  Tensor products of maps are depicted by horizontal juxtaposition.

For every linear map $f\colon V\to V$, we define its complex conjugate $f^*\colon V^*\to V^*$ and transpose $f\Ctrans\colon V^*\to V^*$ by the compositions
$$f\Ctrans=\Diagram{\DMor{square}[-/>] (1,0.5) {$f$}\draw [late arrow] (1,1.5) .. controls +(up:0.5) and +(up:0.5) .. (2,1.5) -- (2,-0.5);\draw [late arrow] (0,1.5) -- (0,-.5) .. controls +(down:0.5) and +(down:0.5) .. (1,-.5);},\qquad f^*=(f^\dag)\Ctrans.$$
Note that the latter can also be equivalently defined by $f^*(e_i^*)=(f(e_i))^*$.

Recall that a {\it compact matrix quantum group} is a pair $G = (A,u)$, where $A$ is a unital C*-algebra, and $u = [u_{ij}] \in M_n(A)$ is a matrix (called a {\it fundamental representation of $G$}) such that (1) the entries $u_{ij}$ of $u$ generate $A$ as a C*-algebra, and (2) both $u$ and its {\it complex conjugate} $\bar u = [u_{ij}^*]$ are invertible in $M_n(A)$.  See \cite{Wor95, NT13, BS09} for more details.  As is standard, we often write $A = C(G)$, regarding $A$ as a possibly non-commutative C$^*$-algebra of continuous functions on the quantum group $G$.  The comultiplication $\Delta:C(G) \to C(G) \otimes C(G)$ is the unital $*$-homomorphism given on generators by $\Delta(u_{ij}) = \sum_{k=1}^n u_{ik} \otimes u_{kj}$. We also denote by $\mathcal O(G) \subseteq C(G)$ the *-algebra generated by the elements $u_{ij}$.  $\mathcal O(G)$ is a Hopf $*$-algebra  with comultiplication $\Delta$ as above, counit $\epsilon:\mathcal O(G) \to \C$ given by $\epsilon(u_{ij}) = \delta_{ij}$, and antipode $S:\mathcal O(G) \to \mathcal O(G)$ given by $(\text{id} \otimes S)(u) = u^{-1}.$             

Every compact matrix quantum group $G$ defines a concrete representation category $\text{Rep}(G)$ by considering representation intertwiners between tensor powers of the fundamental representation $u$ and its complex conjugate $\bar u$ as morphisms. Here, we identify $V$ with the vector space where $u$ acts and $V^*$ with the vector space where $\bar u$ acts. Now, taking two words $w_1$, $w_2$ over two-letter alphabet representing $u$ and $\bar u$, the claim is that $\Cat(w_1,w_2)=\{T\colon V^{\otimes w_1}\to V^{\otimes w_2}\mid Tu^{\otimes w_1}=u^{\otimes w_2}T\}$ defines a concrete representation category. According to Woronowicz--Tannaka--Krein duality \cite{Wor88}, the converse works as well: For every concrete representation category $\Cat$, the $*$-algebra generated by elements $u_{ij}$ subject to relations $Tu^{\otimes w_1}=u^{\otimes w_2}T$, $T\in\Cat(w_1,w_2)$ is actually a Hopf $*$-algebra that defines a compact matrix quantum group with fundamental representation $u=(u_{ij})$. See e.g. \cite{TW17} for more details.

Starting with the work \cite{BS09}, a lot of effort was put into studying quantum groups corresponding to categories based on spiders (also known as categories of partitions).  Examples include the free unitary quantum group $U_n^+$ (corresponding to the smallest category $\langle\rangle=\langle\Rcap,\Lcap\rangle$), the free orthogonal quantum group $O_n^+$ (corresponding to $\langle\Diagram{\Dmor{bcirc}[</<] (0,0.5)}\rangle$, so the arrows can actually be removed obtaining just $\langle \Dcap\rangle$ -- the Temperley--Lieb category), or the free symmetric quantum group $S_n^+$ (corresponding to $\langle\Diagram{\Dmor{bcirc}2/1 (0,0.5)}\rangle$).

\subsection{Complementary spiders}

In this work, we follow on \cite{GHad}, where \emph{complementary spiders} are introduced to the picture. Orthonormal bases $(e_i)$ and $(f_j)$ are called \emph{mutually unbiased} if $|\langle e_i,f_j\rangle|=\frac{1}{\sqrt n}$ for all $i,j$. By the anti-isomorphism $V\to V^*$, this also defines the conjugate basis $(f_j^*)$ of $V^*$. We introduce \emph{white spiders} the same way as we introduced black spiders above replacing black point by white and the $(e_i)$ bases by $(f_j)$. That is,
$$\Diagram{\Dmor{circ}[<.</>.>] (0,.5)}=\sum_i f_i\otimes\cdots\otimes f_i\otimes f_i^\dag\otimes\cdots\otimes f_i^\dag,$$

Let us denote the transition matrix between the two bases by $\frac{1}{\sqrt{n}}H$. That is, $H_{ij}=\sqrt{n}\langle e_i,f_j\rangle$, so $f_j=\sum_i \frac{1}{\sqrt n}H_{ij}e_i$ and $f_j^*=\sum_i \frac{1}{\sqrt{n}}\bar H_{ij}e_i^*$.
We can also introduce the graphical notation
$$H=\Diagram{\Dmor{map}[</>] (1,.5)},\quad
  H^\dag=\Diagram{\Dmor{mapA}[</>] (1,.5)},\quad
  H\Ctrans=\Diagram{\Dmor{mapT}[>/<] (1,.5)},\quad
  H^*=\Diagram{\Dmor{mapC}[>/<] (1,.5)}$$

With this notation, we can write
$$\overbrace{\underbrace{\Diagram{\Dmor{circ}[<>./<>.] (0,.5)}}_{k\times}}^{l\times}=n^{-(k+l)/2}
\Diagram{\Dmor{bcirc}[->./->.] (0,.5)
         \Dmor{mapC}[>/<] (-1,1.5)
         \Dmor{map}[-/>] ( 0,1.5)
         \Dmor{mapA}[</>] (-1,-.5)
         \Dmor{mapT}[>/-] ( 0,-.5)}.$$

The condition of being unbiased implies that $|H_{ij}|=1$. Such a matrix which is unitary up to normalization and whose entries are complex units is by Definition~\ref{D.Hadamard} a Hadamard matrix. These two conditions $H\bullet  H^*=J(=H^*\bullet H)$ and $HH^\dag=nI(=H^\dag H)$ can diagrammatically be expressed as
\begin{equation}
\label{eq.had}
\Diagram{\Dmor{bcirc}2/1 (0.5,1.5)
         \Dmor{map}1/1 (0,0.5)
         \Dmor{mapC}1/1 (1,0.5)
         \Dmor{bcirc}1/2 (.5,-.5)}
=
\Diagram{\Dmor{bcirc}0/1 (0,1)
         \Dmor{bcirc}1/0 (0,0)}
=
\Diagram{\Dmor{bcirc}2/1 (0.5,1.5)
         \Dmor{mapC}1/1 (0,0.5)
         \Dmor{map}1/1 (1,0.5)
         \Dmor{bcirc}1/2 (.5,-.5)},
\qquad
\Diagram{\Dmor{mapA}1/1 (0,1)
         \Dmor{map}1/1 (0,0)}
=
n\;\Did
=
\Diagram{\Dmor{map}1/1 (0,1)
         \Dmor{mapA}1/1 (0,0)},
\end{equation}
where the direction of the arrows does not really matter, so we did not draw them.

Recall that we say that the Hadamard matrix is of order $s$ if all its entries are $s$-th roots of unity, so $H^ {\bullet s}=J$. Diagrammatically, this means that
$$
\Diagram{\Dmor{bcirc}[-.-/-] (1,1.5)
         \Dmor{map}1/1 (0,0.5)
         \Dmor{map}1/1 (2,0.5)
         \Dmor{bcirc}[-/-.-] (1,-.5)}
=
\Diagram{\Dmor{bcirc}0/1 (0,1)
         \Dmor{bcirc}1/0 (0,0)}
$$
where the dots stand for $s$ copies of $H$.

Suppose $H$ is a complex Hadamard matrix. When working with concrete representation categories generated by black and white spiders, we will use the index $H$ (for instance, writing $\langle\Diagram{\Dmor{circ}[<>/><] (0,0.5)},\Diagram{\Dmor{bcirc}[<>/><] (0,0.5)}\rangle_H$) to denote that the black and white spiders are interpreted the above described way, where the two bases are related via the Hadamard matrix $H$.

\subsection{Complexifications of the (quantum) hyperoctahedral group}

In this article we aim to complexify the results of \cite{GHad}, where the second author studied quantum symmetries of (real) Hadamard matrices. These quantum symmetry groups are quantum subgroups of the free hyperoctahedral quantum group $H_n^+$ (introduced by  Bichon \cite{Bic04}). This means that, as an ingredient, we need to consider some suitable complexifications of $H_n^+$.

The classical hyperoctahedral group is defined via the wreath product $H_n=\Z_2\wr S_n$. That is, it can be realized via permutation matrices with $\pm1$ entries. A natural way to complexify these is to replace the $\pm 1$ entries by $s$-th roots of unity for some fixed $s$ or by arbitrary complex units. That is, we define the groups $H_n^s=\Z_s\wr S_n$ for $s=\N\cup\{\infty\}$, where $\Z_\infty=\T$.

The quantum analogues of these groups called \emph{quantum reflection groups} $H_n^{s+}$, $s\in \N\cup\{\infty\}$ can be defined using the {\it free wreath product} construction $\wr_*$, as $H_n^{s+}=\Z_s\wr_{*}S_n^+$. Here $S_n^+$ is the quantum permutation group on $n$ letters \cite{Wa98}.   They first appeared in \cite{Bic04} and were studied more in detail in \cite{BV09}. 

There is another way to quantize the group $H_n^\infty$, first considered by Banica \cite{Ban08}, via the \emph{free complexification} operation on compact quantum groups $G \mapsto G \tilstar\T$.  As is standard in the literature, for $s=2$, we write  $H_n^+ := H_n^{2 +}$ for the free hyperoctahedral quantum group.  By free complexification, one obtains 
a larger (i.e., more free) quantum group $H_n^+\tilstar\T\supsetneq \T\wr_* S_n^+ = H_n^{\infty +}$. See \cite{TW17} for a clearer definition of the free complexification operation.
There is actually a classification result\footnote{For readers who are interested in checking the reference, let us point out a difference in the notation. The work \cite{TW17} uses the black and white point $\bullet$ and $\circ$ to denote the two objects/vector spaces for which we use the arrows $\Fid$ and $\Bid$ instead. On the other hand, they do not explicitly write the spider boxes as there are essentially no other morphisms than these.} \cite{TW17} that says that these are all the ``free'' quantizations of the complex hyperoctahedral groups $H_n^s$ that allow an ``easy'' diagrammatic description. Let us now have a look at the actual definitions and their diagrammatic descriptions.

We will denote here the free complexification of the hyperoctahedral group by $K_n^+:=H_n^+\tilstar\T$. It is defined via the C*-subalgebra $C(K_n^+)$ of the free product $C(H_n^+)*C^*\Z$  generated by $u_{ij}:=v_{ij}z$, where $v= (v_{ij})$ is the fundamental representation of $H_n^+$ and $z$ is the generator of $\C\Z$. Note that $\Z$ can actually be replaced by $\Z_2$ or any cyclic group \cite{TW17,Groglue} in the above definition.  The resulting matrix $u = (u_{ij})$ is a fundamental representation of $K_n^+$, the associated  representation category is generated by the spider $\Diagram{\Dmor{bcirc}[<>/><] (1,.5)}$ (or, equivalently by Frobenius reciprocity, $\Diagram{\Dmor{bcirc}[<><>/] (1,.5)}$). This means that $\mathcal O(K_n^+)$ can equivalently be defined as the $*$-algebra generated by elements $u_{ij}$ subject to relation $uu^\dag=1=u^\dag u$ and $T(u\otimes \bar{u})=(u\otimes \bar{u})T$, where $T=\Diagram{\Dmor{bcirc}[<>/><] (1,.5)}$. Explicitly, the latter relation reads
%\sum_{i=1}^n u_{ij}^*u_{ij}=1=\sum_{j=1}^nu_{ij}^*u_{ij},
\begin{equation}\label{eq.Kn}
u_{ik}u_{il}^*\delta_{ij}=u_{ik}u_{jk}^*\delta_{kl}.    
\end{equation}

The quantum reflection group $H_n^{\infty+}$ is a quantum subgroup of $K_n^+$ (written $H_n^{\infty +}\subset K_n^+$), meaning that the fundamental representation of $H_n^{\infty +}$ has more relations that that of $K_n^+$. The representation category of $H_n^{\infty +}$ is generated by the spider $\Diagram{\Dmor{bcirc}[<</>>] (1,.5)}$.
The corresponding relation reads
\begin{equation}\label{eq.Hninf}
u_{ik}u_{il}\delta_{ij}=u_{ik}u_{jk}\delta_{kl}.    
\end{equation}  Note that we have the following equality in the representation category of $H_n^{\infty +}$:
\[\Diagram{\Dmor{bcirc}[<>/><] (1,.5)}=
 \Diagram{
   \Dmor{bcirc}[<-/>>] (1,.5)
   \Dmor{bcirc}[--/>>] (5,.5)
   \draw[->] (5.5,1) .. controls +(up:0.5) and +(up:0.5) .. (2.5,1) -- (2.5,0) .. controls +(down:0.5) and +(down:0.5) ..  (1.5,0);
   \draw [->] (4.5,1) .. controls +(up:0.2) and +(up:0.2) .. (3.5,1) -- (3.5,0);  
   \draw [late arrow] (1.5,1) .. controls +(up:0.7) and +(up:0.7) .. (6.5,1) -- (6.5,0) .. controls +(down:0.3) and +(down:0.3) .. (5.5,0);
   \draw [mid arrow] (7.5,1) -- (7.5,0) .. controls +(down:0.5) and +(down:0.5) .. (4.5,0);
 }
.\]  In particular, relation~\eqref{eq.Kn} for $C(K_n^+)$ holds in $C(H_n^{\infty +})$, along with  many additional ones.   So indeed $H_n^{\infty +} \subset K_n^+$.  Compare \cite[p. 349]{Bic04}.

Finally, for $s \in \mathbb N$, the representation category of $H_n^{s+}$ is generated by $\Diagram{\Dmor{bcirc}[<<.</] (1,.5)}$ with $s$ inputs, which already generates the two spiders mentioned before, via the equation \[\Diagram{\Dmor{bcirc}[<</>>] (1,.5)}=
 \Diagram{
   \Dmor{bcirc}[<<-.-/] (1,.3)
   \Dmor{bcirc}[/>>-.-] (1,.7)
   \draw [mid arrow] (3,1.2) .. controls +(up:0.2) and +(up:0.2) .. (4,1.2) -- (4,-.2) .. controls +(down:0.2) and +(down:0.2) .. (3,-.2);
   \draw [mid arrow] (1,1.2) .. controls +(up:0.5) and +(up:0.5) .. (5,1.2) -- (5,-.2) .. controls +(down:0.5) and +(down:0.5) .. (1,-.2);
 }.\]
The corresponding relation reads (compare again \cite[p.~350]{Bic04})
$$\sum_{i=1}^nu_{ij_1}\cdots u_{ij_s}=\delta_{j_1,\dots,j_s}.$$
We remind the reader that for $s=2$, we write  $H_n^+ := H_n^{2 +}$.

One can explicitly describe the Hopf *-algebras $\mathcal O(G)$ for $G = H_n^{s+}$ or $K_n^+$ in terms of generators and relations as follows: Let $\A$ be a unital $*$-algebra and $s \in \N$.  A matrix $u = (u_{ij}) \in M_n(\A)$ is called 
\begin{itemize}
    \item a {\it complex quantum reflection matrix} if $u, \bar u$ are unitary and each $u_{ij}$ is a normal partial isometry.
    \item an {\it order $s$-quantum reflection matrix}  if $u$ is a quantum reflection matrix and $u_{ij}^s = u_{ij}u_{ij}^*$ for each $1 \le i,i \le n$.
    \item a {\it free quantum reflection} if $u, \bar u$ are unitary and the entries $u_{ij}$ satisfy relation \eqref{eq.Kn} $u_{ik}u_{il}^*\delta_{ij}=u_{ik}u_{jk}^*\delta_{kl}.$
\end{itemize}
With these notions in mind, we have the following universal unital $\ast$-algebras 
\begin{align*}\mathcal O(K_n^+) &= *\text{-alg}\langle u_{ij}, \ 1 \le i,j \le n: u = (u_{ij}) \text{ is a free quantum reflection} \rangle \\
\mathcal O(H_n^{\infty +}) &= *\text{-alg}\langle u_{ij}, \ 1 \le i,j \le n: u = (u_{ij}) \text{ is a complex quantum reflection} \rangle \\
\mathcal O(H_n^{s +}) &= *\text{-alg}\langle u_{ij}, \ 1 \le i,j \le n: u = (u_{ij}) \text{ is an order-$s$-quantum reflection} \rangle
\end{align*} 
See \cite{BV09, TW17} for details.  By replacing universal unital $\ast$-algebras with universal unital C$^*$-algebras above, one obtains the analogous description of $C(G)$.

So, to summarize, we have quantum groups, with associated function algebras arising as successive quotients of each other,  and their corresponding inclusions of representation categories (below we have $s < \infty$)
$$\begin{array}{rcccl}
H_n^{s+}&\subsetneq&H_n^{\infty+}&\subsetneq&K_n^+,\\
\mathcal O (H_n^{s+})&\leftarrow&\mathcal O(H_n^{\infty+})&\leftarrow &\mathcal O(K_n^+),\\
\langle\underbrace{\Diagram{\Dmor{bcirc}[<<.</] (1,.5)}}_{s\times}\rangle
&\supsetneq&
\langle\Diagram{\Dmor{bcirc}[<</>>] (1,.5)}\rangle
&\supsetneq&
\langle\Diagram{\Dmor{bcirc}[<><>/] (1,.5)}\rangle.
\end{array}$$

As a final remark, let us mention that due to the structure $H_n^{s+}=\Z_s\wr_* S_n^+$, the quantum group $H_n^{s+}$ must act on $sn$ points, that is, $H_n^{s+}\subseteq S_{sn}^+$. This action can be described through \emph{sudoku unitaries}. Let us briefly recall the construction:

\begin{lemma}
\label{L.sudoku}
Let $u$ be the fundamental representation of $H_n^{s+}$. For any $i,j=1,\dots,n$ and $p\in\Z$, define $a_{ij}^{(p)}=\frac{1}{s}\sum_{r=0}^{s-1}\omega^{-rp}u_{ij}^r$, where $\omega$ is a primitive $s$-th root of unity. Then $u_{ij}=\sum_{p=0}^{s-1}\omega^{p}a_{ij}^{(p)}$. In addition, the matrix of the block form
$$p = \begin{pmatrix}
    a^{(0)}&a^{(1)}&\cdots&a^{(s-1)}\\
    a^{(s-1)}&a^{(0)}&\cdots&a^{(s-2)}\\
    \vdots&\vdots&\ddots&\vdots\\
    a^{(1)}&a^{(2)}&\cdots&a^{(0)}
\end{pmatrix}$$
is an $sn \times sn$ quantum permutation matrix, and is a faithful unitary representation of $H_n^{s+}$
\end{lemma}
\begin{proof}
    The fact that $p$ is a quantum permutation as well as the validity of inverse formulas relating the $u_{ij}$ and the $a_{ij}^{(p)}$ follow directly from \cite[Theorem 2.3]{BV09}. The only thing left is showing that the matrix $p$ indeed is a representation of $H_n^{s+}$. This means checking that
    $$\Delta(a_{ij}^p)=\sum_{k=1}^n\sum_{t=0}^sa_{ik}^t\otimes a_{kj}^{p-t}.$$
    This is a straightforward computation, in which we can use the fact that matrices $(u_{ij}^k)$ for any fixed $k$ are representations \cite[Theorem~4.3]{BV09}.
\end{proof}

\subsection{Quantum automorphisms of Hadamard matrices}

We now bring Hadamard matrices into our framework.
For a compact matrix quantum group $G=(C(G),u)$ with $n$-dimensional fundamental representation $u = (u_{ij})$, we will use the shorthand notation $p\in G$ to denote any matrix $p = (p_{ij})\in M_n(\A)$, where $\A$ is some (non-zero) unital C$^*$-algebra such that the map $u_{ij}\mapsto p_{ij}$ extends to a unital $*$-homomorphism $\mathcal O(G)\to\A$.  If we write $p,q\in G$, we assume that both $p,q$ have entries in the same algebra $\A$.

Let $s\in\N \cup \{\infty\}$. We say that two complex Hadamard matrices $H$ and $H'$ are
\begin{itemize}
\item \emph{$s$-quantum equivalent} if there are $p,q\in H_n^{s+}$ such that $H'q=pH$ (or equivalently  $q=H'^{-1}pH$),  
\item \emph{strongly quantum equivalent} if they are $\infty$-quantum equivalent,
\item \emph{weakly quantum equivalent} if there are $p,q\in K_n^+$ such that $H'q=pH$ (or equivalently  $q=H'^{-1}pH$).
\end{itemize}
Note that we have the implications $s$-quantum equivalence $\implies$ strong quantum equivalence $\implies$ weak quantum equivalence.  All of the above quantum equivalences are determined by the choice of quantum group $G =H_n^{s+}, H_n^{\infty +},$ or $K_n^+$, so at times we will also refer to these equivalences simply as $G$-quantum equivalences.

As one can readily see, the above definitions can equivalently formulated as saying $H$ and $H'$ are
\begin{itemize}
\item \emph{$G$-quantum equivalent} if there is $p\in G$ such that $H'^{-1}pH\in G$,
\end{itemize}
The above definitions are representation-theoretic: we are asking for the existence of a non-zero representation $\mathcal O(G)  \to \A$; $u_{ij} \mapsto p_{ij}$, which obeys an additional set of algebraic relations involving $H,H'$.  One can equivalently formulate these definitions in terms of the statement that certain quotient *-algebras of $\mathcal O(G)$ are non-zero.  Indeed, $H, H'$ above are $G$-quantum equivalent if and only if the quotient unital *-algebra \begin{align} \label{eq.isomorphism-algebra}\mathcal A_{G}(H,H') = \mathcal O(G)/\langle  H'^{-1}uH \in G  \rangle \end{align}
is non-zero.\footnote{One direction of this claim is clear.  The other direction relies on the work of Bichon \cite{Bic99} which shows that such algebras are non-zero if and only if they have a non-zero C$^*$-representation.}

When we specialize to the case $H=H'$,  the algebra $\mathcal A_G(H,H)$ is always a non-zero Hopf *-algebra quotient of $\mathcal O(G)$.  We write $\mathcal O(\Aut^+_G(H)) = \mathcal A_G(H,H)$, where $\Aut_G^+(H) = G \cap H G H^{-1} \subseteq G$ is the assoicated {\it $G$-quantum-automorphism group of $H$}.  We have  the chain of quantum subgroups 
$$\Aut^+_{H_n^{s+}}(H)\subseteq \Aut^+_{H_n^{\infty+}}(H) \subseteq \Aut^+_{K_n^{+}}(H),$$ with corresponding inclusions of representation categories $\text{Rep}(\Aut_G^+(H))$:
$$
\langle\underbrace{\Diagram{\Dmor{bcirc}[<<.</] (1,.5)}}_{s\times},\underbrace{\Diagram{\Dmor{circ}[<<.</] (1,.5)}}_{s\times}\rangle_H
\supseteq
\langle\Diagram{\Dmor{bcirc}[<</>>] (1,.5)},\Diagram{\Dmor{circ}[<</>>] (1,.5)}\rangle_H
\supseteq
\langle\Diagram{\Dmor{bcirc}[<><>/] (1,.5)},\Diagram{\Dmor{circ}[<><>/] (1,.5)}\rangle_H
.$$
Note that when we abelianize $\mathcal O(G)$, we recover the usual (coordinate function algebras over the) automorphism groups $\Aut_{H_n^s}(H)$ (when $G = H_n^{s+}$) and $\Aut_{H_n^\infty}(H)$ (when $G = H_n^{\infty +}$ or $K_n^+$). 

As explained in \cite{BCE+, GHad, NT13}, the algebras $\mathcal A_G(H,H')$ appearing in \eqref{eq.isomorphism-algebra}, when non-zero, define {\it $\Aut^+_G(H')-\Aut^+_G(H)$-bi-Galois objects}, or {\it linking algebras}.  As was shown by Schauenberg  \cite{Sch04} (see also \cite{Bic14,NT13}), the existence of existence of bi-Galois objects is equivalent to the existence of a canonical unitary monoidal equivalence of categories $F:\text{Rep}(\Aut^+_G(H)) \to \text{Rep}(\Aut^+_G(H'))$.  In our context, the result of interest is the following:
\begin{theorem} \label{thm.monoidal-equiv}
    Two Hadamard matrices $H$ and $H'$ of size $n$ are $G$-quantum equivalent if and only if the corresponding categories $\text{Rep}(\Aut^+_G(H'))$  and $\text{Rep}(\Aut^+_G(H))$ are unitarily monoidally equivalent via a functor $F$ which maps canonical  generating spiders to canonical  generating spiders. 
\end{theorem}

\begin{proof}[Sketch of Proof] One direction of the proof is easy:  Suppose $\mathcal A_G(H,H') \ne 0$. 
Let $p = (p_{ij})$ be the unitary matrix of generators of $\mathcal A_G(H,H')$.  Similarly, let $u, v$ be the fundamental representations of $\Aut^+_G(H'), \Aut^+_G(H)$, respectively.  Then for any $k, l \in \N$, and any morphism $T \in \text{Mor}_{\Aut^+_G(H')}(u^{\epsilon(1)} \otimes \ldots \otimes u^{\epsilon(k)}, u^{\epsilon'(1)} \otimes \ldots \otimes u^{\epsilon'(l)})$ (here $u^\epsilon, u^{\epsilon'}$ is either $u$ or $\bar u$), it follows from the definitions of $u,v,p$, that there is a unique $F(T) \in \text{Mor}_{\Aut^+_G(H)}(v^{\epsilon(1)} \otimes \ldots \otimes v^{\epsilon(k)}, v^{\epsilon'(1)} \otimes \ldots \otimes v^{\epsilon'(l)})$ such that
\begin{align*}
\big(u^{\epsilon'(1)} \otimes \ldots \otimes u^{\epsilon'(l)}\big)^\dag T\big(u^{\epsilon(1)} \otimes \ldots \otimes u^{\epsilon(k)}\big)
= F(T).
\end{align*}
These equations can easily be seen to uniquely determine the claimed functor $F$.  

The converse direction is highly non-trivial, and relies on a Tannaka-Krein type reconstruction theorem for Galois-objects based on unitary fibre functors.  See \cite[Theorem 2.3.11]{NT13} for details.
\end{proof}

The utility of the above theorem is that it allows us to check whether or not two Hadamard matrices $H,H'$ are $G$-quantum equivalent by directly checking (usually by means of combinatorial and graphical calculus methods) whether or not the fibre functor $F$ exists.  This potentially avoids the seemingly more challenging problem of directly showing that the algebra $\mathcal A_G(H,H') \ne 0$.

\section{(Non-)Examples of Quantum Equivalence} \label{sec:examples}

In \cite{GHad}, the graphical approach to quantum equivalence afforded by Theorem \ref{thm.monoidal-equiv} was used to great effect by the second author.  It was proven there that all {\it real} Hadamard matrices of the same size are mutually 2-quantum equivalent. 

Given the constructions of the previous section, it is natural to ask whether similar statements can be formulated for general complex Hadamard matrices of the same size.  Of course, in the complex case, we have seen that there are multiple choices to be made.  For example,  whether or not to work with fixed order Butson matrices, and which type of quantum symmetries to work with.  We will see now that these choices lead to both positive and negative results generalizing the real case.   

\subsection{A positive result}

We first begin with minimal assumptions on a pair $H,H'$ of complex Hadamard matrices -- only assuming that are both of size $n$.  If at the same time we also allow for the most general notion of quantum equivalence ($K_n^+$-quantum equivalence), we arrive at the following generalization of the result from \cite[Theorem 6.3]{GHad}.    

\begin{theorem}
\label{T.weak}
All complex Hadamard matrices of a fixed size are mutually weakly (i.e., $K_n^+$-) quantum equivalent.
\end{theorem}

The proof of Theorem \ref{T.weak} follows along similar lines as the proof of \cite[Theorem 6.3]{GHad}--we realize the category 
$\text{Rep}(\Aut^+_{K_n^+}(H)) = \langle\Diagram{\Dmor{bcirc}[<><>/] (1,.5)},\Diagram{\Dmor{circ}[<><>/] (1,.5)}\rangle_H$
as an abstract diagrammatic category and show that any closed diagram can be reduced to a number which only depends on $n$, and not the specific choice of Hadamard matrix $H$ of size $n$. This means that all possible fibre functors from this abstract category  (denoted $\langle\Diagram{\Dmor{bcirc}[<><>/] (1,.5)},\Diagram{\Dmor{circ}[<><>/] (1,.5)}\rangle_n$ below) to the concrete categories $\text{Rep}(\Aut^+_{K_n^+}(H))$ have the same kernel \cite[Proposition 1.3]{GHad}. So, in particular, all complex Hadamard matrices must be weakly quantum equivalent.

First, note that the black spiders in the category $\langle\Diagram{\Dmor{bcirc}[<><>/] (1,.5)},\Diagram{\Dmor{circ}[<><>/] (1,.5)}\rangle_H$ satisfy all the standard relations for black spiders in unitary partition categories. The same holds for the white spiders as the basis $(f_j)$ given by normalizing the columns of the Hadamard matrix is orthogonal as well. In particular, we have
$\Diagram{\Dmor{bcirc}[<>/] (1,0.5)}=\Rcap=\Diagram{\Dmor{circ}[<>/] (1,0.5)}$
and
$\Diagram{\Dmor{bcirc}[></] (1,0.5)}=\Lcap=\Diagram{\Dmor{circ}[></] (1,0.5)}$.

In addition observe the following:

\begin{lemma}
In $\langle\Diagram{\Dmor{bcirc}[<><>/] (1,.5)},\Diagram{\Dmor{circ}[<><>/] (1,.5)}\rangle_H$, we have
\begin{equation}
\label{eq.reduction}
\Diagram{\Dmor{bcirc}[->/-.-] (0,1)
         \Dmor{circ}[-.-/>-] (0,0)}
=\frac{1}{n}
\Diagram{\Dmor{bcirc}[/-.-] (0,1)
         \Dmor{circ}[-.-/] (0,0)}
=
\Diagram{\Dmor{bcirc}[>-/-.-] (0,1)
         \Dmor{circ}[-.-/->] (0,0)}
\end{equation}
for any complex Hadamard matrix $H\in M_n(\C)$.
\end{lemma} 
\begin{proof}
This follows immediately from eq.~\eqref{eq.had}
\end{proof}

Now, take any $n\neq 0$ and define $\langle\Diagram{\Dmor{bcirc}[<><>/] (1,.5)},\Diagram{\Dmor{circ}[<><>/] (1,.5)}\rangle_n$ to be the abstract diagrammatic category generated by the given diagrams (and $\Rcap$, $\Lcap$) subject to the standard rules for spiders (with $\Diagram{\draw[mid arrow] (1,0.5) .. controls +(up:0.5) and +(up:0.5) .. (2,0.5) .. controls +(down:0.5) and +(down:0.5) .. (1,0.5);}=n$) and eq.~\eqref{eq.reduction}.

\begin{lemma}
\label{L.alternating}
For every (possibly not reduced) diagram in $\langle\Diagram{\Dmor{bcirc}[<><>/] (1,.5)},\Diagram{\Dmor{circ}[<><>/] (1,.5)}\rangle_n$, the strings going out of each black or white spider always have alternating directions. (In particular, there is always an even number of them.)
\end{lemma}
\begin{proof}
This holds for the generators. It is quite easy to see that the property is preserved under the category operations and applying the reduction rules.
\end{proof}

The following lemma then relates our diagrammatic category $\langle\Diagram{\Dmor{bcirc}[<><>/] (1,.5)},\Diagram{\Dmor{circ}[<><>/] (1,.5)}\rangle_n$ with a similar one without arrows $\langle\Diagram{\Dmor{bcirc}4/0 (1,.5)},\Diagram{\Dmor{circ}4/0 (1,.5)}\rangle_n$ defined in \cite{GHad} (where it was also denoted $\mathsf{NCBipartEven}_n$).

\begin{lemma}
\label{L.corresp}
For every diagram in $$\langle\Diagram{\Dmor{bcirc}[<><>/] (1,.5)},\Diagram{\Dmor{circ}[<><>/] (1,.5)}\rangle_n,$$ we obtain a diagram in $\langle\Diagram{\Dmor{bcirc}4/0 (1,.5)},\Diagram{\Dmor{circ}4/0 (1,.5)}\rangle_n$ by forgetting the arrows. Moreover, the former is reduced if and only if the latter is.
\end{lemma}
\begin{proof}
The first part follows from the fact that the generators for the categories differ only by the arrows. Clearly, if two diagrams with arrows are composable, then they are composable also without arrows. As for the tensor product and involution, the arrows do not play a role. Therefore, all the diagrams that can be generated in the first case, must also be elements of the latter category.

The same holds for the reduction rules, so if a diagram in the first category can be reduced, then it can be reduced also in the second one and the result of the reduction is the same. For the converse, we have to go through all the reduction rules. As for merging the spiders, this does not depend on the direction of the arrows so if it can be done in the second category, it can be done in the first as well. As for the rule \eqref{eq.reduction}, assume a black spider is connected to a white one by a pair of strings in the second category. Then, in the first category, the strings must have opposite arrows by Lemma~\ref{L.alternating}. So, we apply the reduction rule here as well.
\end{proof}

\begin{proof}[Proof of Theorem~\ref{T.weak}]
Take any reduced closed diagram in $\langle\Diagram{\Dmor{bcirc}[<><>/] (1,.5)},\Diagram{\Dmor{circ}[<><>/] (1,.5)}\rangle_n$. Then, by forgetting the arrows, we obtain a closed reduced diagram in $\langle\Diagram{\Dmor{bcirc}4/0 (1,.5)},\Diagram{\Dmor{circ}4/0 (1,.5)}\rangle_n$. But we know that the latter is pure, so the diagram actually has to be empty (just a number).
\end{proof}

\subsection{Some negative results}

Since all Hadamard matrices of order $2$ are $2$-quantum equivalent, one might expect that this might hold for general order $s$.  Namely, that all order $s$-Hadamard matrices of the same size are $s$-quantum equivalent. We show below that this is not the case.  We actually conjecture that for every order $s>2$ there exists a pair of complex Hadamard matrices of certain size $n$ and order $s$ that are not $s$-quantum  equivalent.  We show a couple of concrete examples below.

Recall that the Fourier matrix $\F_n\in M_n(\C)$ with entries $[\F_n]_{ij}=\omega^{ij}$, where $\omega$ is a primitive $n$-th root of unity, (the Fourier transform on $\Z_n$) is a complex Hadamard matrix of order $n$.

\begin{proposition}
\label{P.Fourier}
Consider $n=mp^2$ for $m,p\in\N$, $p\neq 1$. Then the order $n$ Hadamard matrices $\F_n$ and $\F_{mp}\otimes\F_p$ are not strongly quantum equivalent. Hence, they are also not $n$-quantum equivalent.
\end{proposition}
\begin{proof}
Recall the categorical characterization of quantum isomorphism from Theorem \ref{thm.monoidal-equiv}: $H$ is strongly quantum equivalent to $H'$ if and only if there is monoidal equivalence
$$\langle\Diagram{\Dmor{bcirc}[<</>>] (1,.5)},\Diagram{\Dmor{circ}[<</>>] (1,.5)}\rangle_H
\to
\langle\Diagram{\Dmor{bcirc}[<</>>] (1,.5)},\Diagram{\Dmor{circ}[<</>>] (1,.5)}\rangle_{H'}$$
mapping generators to generators. To show that $H$ is not quantum equivalent to $H'$, it is enough to find a closed diagram in the above category that evaluates differently when interpreting by $H$ or $H'$.

We can consider the following one, where all the dots stand for $p$ connections in total:
$$\Diagram{
\Dmor{bcirc}[-0->.>/] ( 0,1.5)
\Dmor{circ}[-0-/>.>] (-1.5,0.5)
\Dmor{circ}[>.>/-0-] ( 1.5,0.5)
\Dmor{bcirc}[/>.>-0-] ( 0,-.5)
}=
\Diagram{
\Dmor{bcirc}[-0-/]   (0,2.5)
\Dmor{circ}[-0-/>.>]  (0,1.5)
\Dmor{bcirc}[-0-/>.>] (0,0.5)
\Dmor{circ}[-0-/>.>]  (0,-.5)
\Dmor{bcirc}[/>.>]  (0,-1.5)
\draw [mid arrow] (0,2.5) .. controls +(up:.5) and +(up:.5) .. (2,2.5) -- (2,-1.5) .. controls +(down:0.5) and +(down:0.5) .. (0,-1.5);
}=n^{-2p}
\Diagram{
\Dmor{bcirc}[-.-/]   (0,4.5)
\Dmor{map}[-/>]  (-1,3.5)
\Dmor{map}[-/>]  (1,3.5)
\Dmor{bcirc}[-.-/>.>]  (0,2.5)
\Dmor{mapA}[-/>]  (-1,1.5)
\Dmor{mapA}[-/>]  (1,1.5)
\Dmor{bcirc}[-.-/>.>] (0,.5)
\Dmor{map}[-/>]  (-1,-.5)
\Dmor{map}[-/>]  (1,-.5)
\Dmor{bcirc}[-.-/>.>]  (0,-1.5)
\Dmor{mapA}[-/>]  (-1,-2.5)
\Dmor{mapA}[-/>]  (1,-2.5)
\Dmor{bcirc}[/>.>]  (0,-3.5)
\draw [mid arrow] (0,4.5) .. controls +(up:.5) and +(up:.5) .. (2,4.5) -- (2,-3.5) .. controls +(down:0.5) and +(down:0.5) .. (0,-3.5);
}
=n^{-2p}\Tr((H^{\bullet p}(H^\dag)^{\bullet p})^2)
$$

In order to see that the above element is indeed in the category, the reader should convince themselves that the spider $\Diagram{\Dmor{bcirc}[<</>>] (1,.5)}$ (or $\Diagram{\Dmor{circ}[<</>>] (1,.5)}$) already generates any black (or white) spider with an equal amount of ingoing and outgoing arrows. (See also \cite[Lemma~1.3(c)]{TW18}.)

Now, let us compute the quantity for the two complex Hadamard matrices mentioned in the formulation. For $H=\F_n$, we can easily derive $H^{\bullet p}=\F_{mp}\otimes J_p$, so $H^{\bullet p}(H^\dag)^{\bullet p}=mp^2I_{mp}\otimes J_p$. Squaring it, we get $m^2p^5I_{mp}\otimes J_p$, which has trace $m^3p^7=n^3p$. In contrast, taking $H=\F_m\otimes\F_p$, we get $H^{\bullet n}=\F_m\otimes J_{p^2}$, so $H^{\bullet n}(H^\dag)^{\bullet n}=mp^2 I_m\otimes J_{p^2}$. Squaring it, we get $m^2p^6 I_m\otimes J_{p^2}$, which has trace $m^3p^8=n^3p^2$.
\end{proof}

For prime orders, the situation is more difficult. First, because constructing pairs of classically non-equivalent Hadamard matrices is more involved (there is just one Fourier matrix of order $p$, so we need to construct something else). Secondly, one has more reduction rules: If $H$ is a Hadamard matrix of a prime order $p$, then $H^{\bullet k}$ is also a Hadamard matrix of order $p$ for any $k\perp p$. Hence, not only do we have $\Diagram{\Dmor{circ}[</>] (1,0.5)}=\Fid$, but also
$\Diagram{\Dmor{bcirc}[-.-/-.-] (1,1.5)
          \Dmor{circ}[-.-/>.>] (1,.5)
          \Dmor{bcirc}[-.-/>.>] (1,-.5)
}=\Diagram{\Dmor{bcirc}[-.-/-.-] (1,.5)}$. So, for instance, the diagram from the proof above can always be reduced.

Nevertheless, as an evidence that also for prime orders one can have strongly quantum non-isomorphic Hadamard matrices, we mention the following.

\begin{proposition}
\label{P.H7}
There are Hadamard matrices of order 7 and size 14 that are not strongly quantum equivalent.
\end{proposition}
\begin{proof}
Consider the diagram
$$
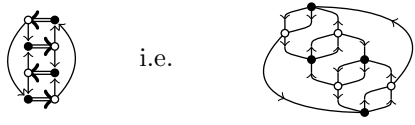
\Diagram{
\draw[->,shorten >=2pt] (0,2) .. controls (-1,1) and (-1,0) .. (0,-1);
\draw[->,shorten >=2pt] (1,-1) .. controls (2,0) and (2,1) .. (1,2);
\draw[->,shorten >=2pt] (0,0) -- (0,-1);
\draw[->,shorten >=2pt] (1,-1) -- (1,0);
\draw[->,shorten >=2pt] (0,0) -- (0,1);
\draw[->,shorten >=2pt] (1,1) -- (1,0);
\draw[->,shorten >=2pt] (0,2) -- (0,1);
\draw[->,shorten >=2pt] (1,1) -- (1,2);
\draw [double,->] (0,-1) node[bcirc]{} -- (1,-1) node[circ]{};
\draw [double,->] (1, 0) node[bcirc]{} -- (0, 0) node[circ]{};
\draw [double,->] (0, 1) node[bcirc]{} -- (1, 1) node[circ]{};
\draw [double,->] (1, 2) node[bcirc]{} -- (0, 2) node[circ]{};
}
\qquad\hbox{i.e.}\qquad
\Diagram{
	\draw[mid arrow] (-2,1.5) .. controls +(-2,-2) and +(-3,0) .. (1,-1.5);
	\draw[mid arrow] (2,-.5) .. controls +(2,2) and +(3,0) .. (-1,2.5);
	\Dmor{bcirc}[/->>]  (1,-1.5)
	\Dmor{circ}[>/>--]  (0,-.5)
	\Dmor{circ}[--0/>]  (2,-.5)
	\Dmor{bcirc}[-/->>] (-1,.5)
	\Dmor{bcirc}[>>-/-] (1,.5)
	\Dmor{circ}[>/0--]   (-2,1.5)
	\Dmor{circ}[-->/>] (0,1.5)
	\Dmor{bcirc}[>>-/]  (-1,2.5)
}
$$

There are three classically nonequivalent complex Hadamard matrices of order seven and size fourteen \cite{LOS20}. If we compute the invariant corresponding to the above given diagram times $2^57^3$ for them, we get the values $3$, $23$, and $3$, so the second one is not strongly quantum equivalent to the other two.
\end{proof}

\begin{remark}
    Butson in his original work~\cite{But62} figured out a way how to construct a Hadamard matrix of order $p$ and size $2p$. It starts with a Fourier matrix $\F_p$ and somehow doubles its size while keeping the order $p$. See also \cite{Szo14} for a much clearer description of the construction.

    This construction depends on an additional number $r$, which is supposed to be a quadratic nonresidue mod $p$. (The construction works in general, but the result is of order $2p$ if $r$ is a quadratic residue.) The construction may or may not produce classically equivalent Hadamard matrices for different $r$. For $p=5$, the two results actually are classically equivalent. But in general, it appears that the results are typically not strongly quantum equivalent. To be more precise, we checked by direct computation (using a computer) that for all prime numbers $p=7,11,\dots,47$, this construction gives examples of complex Hadamard matrices of order $p$ and size $2p$ that can have different values of the invariant mentioned in the proof of Proposition~\ref{P.H7} and hence are not mutually strongly quantum equivalent. However, we were not able to prove it for general $p$. Note also an interesting fact that the invariant times $2^5p^3$ appears to almost always be a natural number except for a few exceptions when it is complex.
\end{remark}

\section{Hadamard equivalence game}
\label{sec:games}

We now study quantum symmetries of Hadamard matrices through the lens of non-local games.  More precisely, we define a synchronous non-local game whose local strategies captures the $s$-equivalence of a pair of Hadamard matrices.  Quantum commuting strategies for this game will capture the notion of quantum $s$-equivalence.  We assume the reader is familiar with the basic concepts of (synchronous) non-local games, their strategy classes, and algebras associated to synchronous games.  For more details, we refer the reader to \cite{HMPS19} and the excellent survey article \cite{HP25}.   

\subsection{The game}

 Suppose $H$ and $H'$ are complex Hadamard matrices of order $s\in\N$ and size $n$. To simplify the notation, let us denote $[n]=\{1,\dots,n\}$ and recall that $\hat\Z_s$ denotes the group of all $s$-th roots of unity. In the game described below, Alice and Bob are spatially separated (non-communicating) players who  are cooperating to try and convince the verifier (referee) that there is a (quantum) $s$-equivalence $H$ to $H'$. Let us denote this game by $\mathcal G(H,H') = (I,O, \lambda)$. Here $I$ denotes this inputs (game questions), $O$ denotes the outputs (game answers), and $\lambda:O \times O \times I \times I \to \{0,1\}$ denotes the predicate (rule function of the game).     

More precisely, one round of $\mathcal G(H,H')$ is played as follows: Each player receives a question from the input set $I=\{r_i,c_i\}_{i=1}^n$, which should be interpreted as a specification of the $i$-th row or column of $H$. Their answers are from the output set $O=I\times\hat\Z_s$. The idea is that they should answer where does the given equivalence map the given row/column of $H$ to,  and by which complex unit is it supposed to be multiplied. The rules for winning this round of the game are as follows
\begin{itemize}
    \item ({\it Synchronicity}) If the two players get exactly the same question, they have to give exactly the same answer. 
    \item If the question is a row (resp. column), then the answer must be a row (resp. column).
    \item If one player gets $r_i$ as a question and answers $(r_k,a)$ and if the other gets $c_j$ as question and answers $(c_l,b)$, they win if and only if 
    $$aH_{ij}=bH'_{kl}.$$
\end{itemize}

Note that we do not give any rule on what to answer if the players are asked two different rows (or columns). In this case the players automatically win the round on any such pair of inputs.  

Our main result is the following.

\begin{theorem}
\label{T.HadGame}
    For the $s$-equivalence game $\mathcal G(H,H')$, we have \begin{itemize}
        \item $\mathcal G(H,H')$ has a perfect local (loc) strategy if and only if $H,H'$ are classically $s$-equivalent.
        \item  $\mathcal G(H,H')$ has a perfect quantum commuting (qc) strategy if and only if $H,H'$ are quantum $s$-equivalent.
        \item $\mathcal G(H,H')$ has a perfect quantum (q) strategy if and only if the linking algebra $\mathcal A_{H_n^{s+}}(H,H')$ has a finite dimensional representation.
        \item $\mathcal G(H,H')$ has a perfect quantum approximate(q) strategy if and only if the linking algebra $\mathcal A_{H_n^{s+}}(H,H')$ has a non-zero homomorphism into $R^\omega$ (the ultrapower of the hyperfinite II$_1$-factor).
    \end{itemize}
\end{theorem}

We will prove this Theorem by establishing a $*$-isomorphism between the synchronous game algebra $\mathcal A(\mathcal G(H,H'))$ and the linking algebra 
$\mathcal A_{H_n^{s+}}(H,H')$.  In fact, we will prove this isomorphism in a more general framework based the notion of monomial isomorphism games.

\subsection{Monomial isomorphisms}

Consider a pair of matrices $A,B\in M_n(\C)$. Consider a set $M\subseteq\C$. We say that $A$ is $M$-monomially $t$-isomorphic ($t={\rm loc,q,qc,\dots}$) to $B$ if there is a perfect $t$-strategy to the following (bisynchronous) game $\mathcal G = (I,O, \lambda)$: The input and output sets are given by $I=[n]:=\{1,\dots,n\}$, $O=[n]\times M$. The rule function $\lambda: I \times I \times O \times O \to \{0,1\}$ for winning a round of the game is given by 
\begin{align*}
\lambda(i,j\mid k,a,l,b)&=\delta_{aA_{ij},bB_{kl}}\delta_{k\neq l}\qquad\text{for $i\neq j$,}\\
\lambda(i,i\mid k,a,l,b)&=\delta_{aA_{ii},aB_{kk}}\delta_{k,l}\delta_{a,b}.
\end{align*}

We will focus here on the case $M=\hat\Z_s$; we then refer to $s$-monomial $t$-isomorphism instead of $\hat\Z_s$-monomial $t$-isomorphism.

In addition, we associate to any matrix $A \in M_n(\C)$ its {\it $s$-monomial quantum automorphism group}.  This is the quantum subgroup of $H_n^{s+}$ given by adding the relation $uA=Au$, where $u$ is the fundamental representation of $H_n^{s+}$. This naturally induces a notion of a quantum isomorphism. To be more concrete, we say that $A$ is $s$-monomially quantum isomorphic to another matrix $B$ if there is a non-zero *-algebra generated by variables $u_{ij}$, $1 \le  i,j \le n$ satisfying
\begin{itemize}
\item $u$ is unitary,
\item $u_{ik}u_{jk}^*\delta_{kl}=u_{ik}u_{il}^*\delta_{ij}$,
\item $\sum_j u_{i_1j}\cdots u_{i_sj}=\delta_{i_1,\dots,i_s}$,
\item $uA=Bu$.
\end{itemize}
As in the case of quantum equivalences of Hadamard matrices, we define the universal $*$-algebra with generators $u_{ij}$ and relations as above to be the \emph{$s$-monomial quantum isomorphism algebra}.

\begin{theorem}
    Matrices $A$ and $B$ are $s$-monomially quantum isomorphic (by the algebraic definition) if and only if they are $s$-monomially qc-isomorphic (by the non-local game definition). 
\end{theorem}
\begin{proof}
By \cite{HMPS19}, it is enough to study the $*$-algebra of the game $\mathcal G$ and show that it is isomorphic to the $s$-monomial quantum isomorphism algebra. The algebra of the game is defined to be the universal $*$-algebra $\mathcal A(\mathcal G)$ given by generators $e_{kai}$  subject to relations
\begin{itemize}
    \item $e_{kai}^*=e_{kai}=e_{kai}^2$,
    \item $\sum_{k,a}e_{kai}=I$,
    \item $e_{kai}e_{lbj}=0$ if $\lambda(i,j\mid k,a,l,b)=0$.
\end{itemize}

We start by showing  the map 
$$u_{ki}\mapsto e_{ki}:=\sum_{a\in\hat\Z_s}ae_{kai}$$
extends to a $*$-homomorphism from the $s$-monomial quantum isomorphism algebra $\mathcal A(\mathcal G)$. In order to do that, we have to show that the $e_{ki}$ satisfy the defining relations of the $u_{ki}$.

Note first that the game rules feature a slightly stronger version of bisynchronicity -- having just $k=l$ already forces $a=b$ and $i=j$.  So, we have $e_{kai}e_{kbj}=\delta_{ij}\delta_{ab}e_{kai}$. Using this, we can use similar argumentation as in \cite[Section~2]{PR21} to obtain $$\sum_{i,a}e_{kai}=1.$$

Now, we are ready to prove all the relations. In the following, we denote by $e=(e_{ki})$ the matrix made out of the elements $e_{ki}$.
We have that
\begin{align*}
[ee^\dag]_{kl}&=\sum_{i}e_{ki}e_{li}^*=\sum_{i,a,b}a\bar b e_{kai}e_{lbi}=\sum_{i,a}e_{kai}\delta_{kl}=\delta_{kl},\\
[e^\dag e]_{ij}&=\sum_{k}e_{ki}^*e_{kj}=\sum_{k,a,b}\bar abe_{kai}e_{kbj}=\sum_{k,a}e_{kai}\delta_{ij}=\delta_{ij},\\
e_{ki}e_{kj}^*&=\sum_{a,b}a\bar be_{kai}e_{kbj}=\sum_a e_{kai}\delta_{ij}=e_{ki}\delta_{ij},\\
e_{ki}e_{li}^*&=\sum_{a,b}a\bar be_{kai}e_{lbi}=\sum_a e_{kai}\delta_{kl}=e_{ki}\delta_{kl},\\
\sum_j e_{i_1j}\cdots e_{i_sj}&=\sum_{j,a_1,\dots,a_s}a_1\cdots a_se_{i_1a_1j}\cdots e_{i_sa_sj}\\&=\sum_{j,a}e_{i_1aj}\delta_{i_1,\dots,i_s}=\delta_{i_1,\dots,i_s},\\
[eA]_{kj}&=\sum_i e_{ki}A_{kj}=\sum_{a,i}ae_{kai}A_{ij}\underbrace{\sum_{l,b}e_{lbj}}_{=I}=\sum_{a,b,i,l}aA_{ij}e_{kai}e_{lbj}\\&=\sum_{a,b,i,l}bB_{kl}e_{kai}e_{lbj}=\sum_{i,a}e_{kai}\sum_{l,b}bB_{kl}e_{lbj}=\sum_l B_{kl}e_{lj}\\&=[Be]_{kj}.
\end{align*}

Now, let us consider the converse direction.  That is, we construct the corresponding inverse $*$-homomorphism from  $\mathcal A(\mathcal G)$ to the $s$-monomial quantum isomorphism algebra.   Since the generators $u_{ki}$ of the $s$-monomial quantum isomorphism algebra satisfy the relations of the generators of $\mathcal O( H_n^{s+})$, we can use Lemma~\ref{L.sudoku} to find projections $u_{kai}$, $1 \le i,k \le n, a \in \hat \Z_s$, such that $u_{ki}=\sum_a au_{kai}$, $\sum_{k,a}u_{kai}=1$, and $\sum_{i,a}u_{kai}=1$. It is now enough to check that we have $u_{kai}u_{lbj}=0$ whenever $\lambda(i,j\mid k,a,l,b)=0$. We now go through the cases:  If $i\neq j$ and $k=l$, then $u_{ki}u_{lj}=0$ because the $s$-monomial quantum isomorphism algebra is a quotient of $\mathcal O(H_n^{s+})$.  In particular, this gives
\[0 = u_{ki}u_{lj}=\sum_{a',b' \in \hat \Z_s}a'b'u_{ka'i}u_{kb'j}.\]
Multiplying this equation on the left by a fixed $u_{kai}$ and on the right by a fixed $u_{kbj}$ and using the fact that  $\{u_{kai}\}_{k,a}$ and $\{u_{kai}\}_{i,a}$ are projection valued measures (hence orthogonal), we obtain $abu_{kai}u_{kbj} = 0 \implies u_{kai}u_{kbj} = 0$. The same arguments apply in the case $i=j$ and $k\neq l$. It remains to show that $u_{kai}u_{lbj}=0$ whenever $aA_{ij}\neq bB_{kl}$. This relation is a consequence of the following equality:
$$aA_{ij}u_{kai}u_{lbj}=\sum_{c,x}u_{kai}cA_{xj}u_{kcx}u_{lbj}=\sum_{d,y}u_{kai}dB_{ky}u_{ydj}u_{lbj}=bB_{kl}u_{kai}u_{lbj},$$
where we used that $uA=Bu$, which means $\sum_{c,x}cu_{kcx}A_{xj}=\sum_{d,y}B_{ky}du_{ydj}$.  

To sum up, the above relations on the $u_{kai}$ ensure that the map $e_{kai}\mapsto u_{kai}$ extends to a unital  $*$-homomorphism from $\mathcal A(\mathcal G)$ to the $s$-monomial quantum isomorphism algebra.  This is clearly an inverse to the previously defined morphism going the other direction, and hence $\mathcal A(\mathcal G)$ and the  $s$-monomial quantum isomorphism algebra are $*$-isomorphic.   
\end{proof}

Focusing on Hadamard matrices provides a proof of Theorem~\ref{T.HadGame}:

\begin{proof}[Proof of Thm.~\ref{T.HadGame}]
    Observe that there is a difference between the quantum equivalence notions for Hadamard matrices of the previous section and those that were introduced in this section. Namely that here we assume that there is a single signed quantum permutation matrix acting on a given  matrix $A$ from both sides. There is an analogous difference in the definition of the game. This can be overcome as follows.

    For a Hadamard matrix $H$, we define a matrix
    $$A:=\begin{pmatrix}I&H\\H^\dag&0\end{pmatrix}.$$
    Then the quantum monomial $s$-isomorphism defined above translates to the quantum $s$-equivalence of complex Hadamard matrices. To see that, first note that if $A_{ii}\neq B_{jj}$, then a quantum isomorphism $A\to B$ must satisfy $u_{ij}=0$ (considering $\sum_i u_{ki}A_{ij}=\sum_l B_{kl}u_{lj}$, multiply it by $u_{kj}^*$ and get $u_{kj}u_{kj}^*A_{jj}=u_{kj}u_{kj}^*B_{kk}$ -- this is a standard argument from the theory of graph quantum symmetries, cf. \cite[Lemma~3.2.3]{Ful06}). Consequently, taking two Hadamard matrices $H_1$ and $H_2$ and the corresponding matrices $A$, $B$ constructed as above, the corresponding quantum isomorphism must be of the form $u=\begin{pmatrix}p&0\\ 0&q\end{pmatrix}$. But now the relation $uA=Bu$ just means that $pH=Hq$ (and $qH^\dag=H^\dag p$ which is equivalent to that).

    As for the games, we have to prove that the Hadamard equivalence game for Hadamard matrices $H$ and $H'$ has a perfect strategy if and only if the monomial equivalence game for the corresponding matrices $A$ and $B$ has a perfect strategy. We will do that by showing that the ``relevant'' rules of the two games coincide. By this, we mean the following: Observe that if $H$ is $n\times n$, then $A$ is $2n\times 2n$ and the indices for the latter can be labelled as $\{r_i,c_i\}_{i=1}^n$. Hence, in the two cases we have the same question and answer sets. Now if we are in the particular situation, when one of the players gets a row question and gives a row answer and the other player gets a column question and gives a column answer, then the two games coincide. This is because $A_{r_ic_j}=H_{ij}=A_{c_jr_i}$, so
    $$aA_{r_ic_j}=bB_{r_kc_l}\quad\Leftrightarrow\quad aH_{ij}=bH'_{kl}\quad\Leftrightarrow\quad aA_{c_jr_i}=bB_{c_lr_k}.$$

    And we claim that this is enough because in both games every perfect strategy has to give a row answer on a row question and a column answer on a column question. Indeed, for the latter game this is proven by the block diagonal structure of the quantum isomorphisms we showed above; as for the former one, it is a consequence of the first game rule (see Lemma~\ref{L.ns}). For the case that both players get a row question (or a column question) and give a row answer (resp. column answer), there are no rules in either game except for the synchronicity.
\end{proof}

\begin{question}
We focused on the case $M=\hat\Z_s$, but in principle, the set $M\subseteq\C$ could be arbitrary. In particular, $M=\T$ should correspond to the notion of strong quantum equivalence. Non-local games with uncountable output set have only recently been introduced and studied in \cite{BTT24}. In particular, at this time there seems to be no notion of a game algebra in this setting, so proving the above theorems for this case is not straightforward. We leave this as an open problem.
\end{question}

\subsection{A non-synchronous game}
There is also a natural non-synchronous version of the above game whose perfect deterministic strategies correspond to monomial $s$-equivalence between complex Hadamard matrices. We briefly outline this game below. We were, however, unable to prove that it is equivalent to the synchronous version.  That is, we do not know if perfect quantum commuting strategies for this non-synchronous game correspond to our algebraic notion of quantum $s$-equivalence of Hadamard matrices. 

The rules of the game are as follows: The verifier sends an index $i\in[n]$ to Alice and $j\in[n]$ to Bob. Alice answers a pair $(k,a)$ with $k\in[n]$ and $b\in\hat\Z_s$. Likewise, Bob has to answer a pair $(l,b)\in[n]\times\hat\Z_s$. They win the round of the game if $aH_{ij}=bH'_{kl}$, and loose otherwise.

\begin{proposition}
    Perfect deterministic strategies in this game are exactly the monomial $s$-equivalences between $H$ and $H'$.
\end{proposition}
\begin{proof}
    First, observe that a deterministic strategy is given by maps $\pi,\sigma\colon[n] \to [n]$ and $a,b\colon [n] \to\hat \Z_s$ such that $a(i)H_{ij}=b(j)H'_{\pi(i)\sigma(j)}$.
    
    On the other hand, recall that a monomial equivalence is given by a pair of monomial matrices $P,Q$ with entries in $\hat\Z_s$ such that $H=P^{-1}H'Q$. Suppose $Pe_i=a_ie_{\pi(i)}$ and $Qe_i=b_ie_{\sigma(i)}$ with $a_i,b_i\in\hat\Z_s$ and $\pi,\sigma\in S_n$. Then the relation $H=P^{-1}H'Q$ reads $H_{ij}=a_i^{-1}b_jH'_{\pi(i)\sigma(j)}$, which is exactly the formula above.

    The only subtlety that remains to explain is that we assumed $\pi,\sigma$ are bijections in the latter case, but not in the former. But since both $H,H'$ are scalar multiples of unitaries, they cannot be related by non-bijective map in the above formula.
\end{proof}

\begin{question}
    Is it possible to algebraically characterize perfect quantum strategies for this non-synchronous game? Does it hold that there is a perfect $t$-strategy for this game if and only if there is a perfect $t$-strategy for the synchronous version?
\end{question}

We show below that at least the right-left implication holds.

\subsection{Remarks on synchronizing a game}

Suppose we have a non-local game $\mathcal G$ with input sets $I_1$, $I_2$ and output sets $O_1$, $O_2$. Denote by $\lambda\colon I_1\times I_2\times O_1\times O_2\to\{0,1\}$ the boolean predicate characterizing conditions for winning the game $\mathcal G$. Below we construct a synchronized version $\mathcal G^s = (I,O,\lambda\sync)$, of this game.  Based on conversations with experts, the following construction seems likely to be well known, but we could not find a reference in the literature.  In particular, we did not find any answers to Question \ref{quest:opposite} below.

We define input and output sets as disjoint unions $I=I_1\sqcup I_2$, $O=O_1\sqcup O_2$. For an element $i\in I_\alpha$, $\alpha=1,2$, we denote by $i^\alpha$ the corresponding element of $I$. Likewise for $O$. Now, put
\begin{align*}
    \lambda\sync(i^\alpha,i^\alpha\mid k^\beta,l^\gamma)&=\delta_{kl}\delta_{\alpha\beta\gamma}\\
    \lambda\sync(i^\alpha,j^\alpha\mid k^\beta,l^\gamma)&=1\qquad\text{assuming $i\neq j$}\\
    \lambda\sync(i^1,j^2\mid k^1,l^2)&=\lambda(i,j\mid k,l)=\lambda\sync(j^2,i^1\mid l^2,k^1)\\
    \lambda\sync(i^\alpha,j^\beta\mid k^\beta,l^\alpha)&\qquad\text{for $\alpha\neq\beta$ can be chosen arbitrarily}
\end{align*}
We call $\mathcal G^s$ the synchronous version of $\mathcal G$.

\begin{remark}
    The synchronous version of the Hadamard $s$-equivalence game is an instance of the synchronization procedure applied to the non-synchronous version.
\end{remark}

The following lemma explains, why the concrete choice for the last rule is irrelevant.

\begin{lemma}
\label{L.ns}
For any perfect non-signalling strategy for the game $\mathcal G^s$, it holds that, given any question $i\in I_\alpha$, to Alice (or Bob) she (he) will give an answer from $O_\alpha$ with probability 1.
\end{lemma}
\begin{proof}
    By definition of non-signalling strategies, the Alice's marginal probability is independent of Bob's input, so assuming $\beta\neq\alpha$, we have $p_A(k^\beta\mid i^\alpha)=\sum_{l^\gamma}p_A(k^\beta,l^\gamma\mid i^\alpha,i^\alpha)=0$.
\end{proof}

\begin{remark}
    If the strategy is quantum and $(H,\psi,(E_{k^\beta i^\alpha}),(F_{l^\beta,j^\alpha})$ is it's minimal realization, then the lemma implies that $E_{k^\beta i^\alpha}=0$ for $\alpha\neq\beta$, which can also be deduced directly from the rule $\lambda\sync(i^\alpha,i^\alpha\mid k^\beta,k^\beta)=0$
\end{remark}

\begin{proposition}
    There is a perfect $t$-strategy for the original game $\mathcal G$ if there is a perfect $t$-strategy for the synchronous version $\mathcal G^s$, for $t={\rm loc, q, qc}$.
\end{proposition}
\begin{proof}
    Suppose there is a perfect quantum strategy for $\lambda\sync$ and $(H,\psi,(E\sync_{k^\beta i^\alpha}),(F\sync_{l^\beta,j^\alpha})$ its minimal realization. (Assuming $H$ is finite-dimensional in the $t=\rm q$ case and all measurements mutually commuting in the $t=\rm loc$ case.) We claim that $H$, $\psi$ together with operators $E_{ki}=E_{k^1i^1}\sync$, $F_{lj}=F_{l^2j^2}\sync$ provide a quantum strategy for $\lambda$. We have $\sum_k E_{ki}=I$ thanks to the fact that $E_{l^2i^1}=0$ by the above lemma/remark. Same for $F$. Having $0=\lambda(i,j\mid k,l)=\lambda\sync(i^1,j^2\mid k^1,l^2)$, we must have $0=\langle\psi,E_{i^1k^1}\sync F_{j^2l^2}\sync\psi\rangle=\langle\psi,E_{ik}F_{jl}\psi\rangle$ which is all we need.
\end{proof}
We close this section with an open problem related to the above construction.
\begin{question} \label{quest:opposite}
    Does the opposite implication hold in the above proposition? If not, find a counterexample.
\end{question}

\section{Hadamard graphs} \label{sec:graphs}

In the real case, given a Hadamard matrix $H$ of size $n$, one can construct its \emph{Hadamard graph} on $4n$ vertices, whose symmetries exactly correspond to the symmetries of the original Hadamard matrix \cite{McK79}. As shown by the second author in \cite{GHad}, the same holds for quantum symmetries.

In the complex (Butson) case, it turns out that the construction of the associated Hadamard graph can be generalized, as we explain below.  As far as we know, the correspondence between Butson Hadamard matrices and their associated graphs was not formally studied anywhere in the literature. The reason for this might be that the correspondence between these objects at the level of (quantum) symmetries does not work so nicely anymore. As we show below, only one inclusion still works in general, and this inclusion holds for both classical and quantum symmetries.

Starting with a monomial matrix $M$ of size $n$ with entries in $\hat\Z_s=\{\omega^i\}_{i=0}^{s-1}$, $\omega=e^{2\pi i/s}$, we construct a $sn \times sn$ matrix $\pi(M)$ by replacing each $\omega^k$ with an $s\times s$ cyclic permutation matrix $C^k$ (where $C$ corresponds to the cycle $(1\;2\;\cdots\;s)$). In particular, starting with a complex Hadamard matrix $H$ of order $s$, we obtain an incidence matrix $D=\pi(H)$ of a symmetric class regular transversal design over $\Z_s$. The association $H \mapsto \pi(H)$ is actually a two-way correspondence that works more generally with generalized Hadamard matrices with group entries. See \cite[Section~VIII.3]{BJL} for details.

Without going into the intricacies of design theory, let us just mention that designs are combinatorial incidence structures that can be characterized by their incidence matrices, which are rectangular matrices consisting of zeros and ones and should satisfy some additional requirements. Two designs with incidence matrices $D_1$ and $D_2$ are said to be isomorphic if there are permutation matrices $\hat P$ and $\hat Q$ such that $D_1=\hat PD_2\hat Q^{-1}$. Considering the construction from above $D_i=\pi(H_i)$, we see that equivalences of complex Hadamard matrices induce isomorphisms of the corresponding designs as we can take $\hat P=\pi(P)$, $\hat Q=\pi(Q)$ assuming that $H_1=PH_2Q^{-1}$, $P,Q\in H_n^s$.

Unfortunately, the converse is not true as it is not always guaranteed that $\hat P$ and $\hat Q$ are of the specific form $\pi(P)$, resp. $\pi(Q)$. For instance, take any complex Hadamard $H$ with $2<s<\infty$, $D:=\pi(H)$. Now, for any automorphism $\beta\colon\Z_s\to\Z_s$, we can consider its permutation matrix $B$ and get $BC^kB^{-1}=C^{\beta(k)}$. Put $\hat P=I_n\otimes B$, $\hat Q=I_n\otimes B$, $D'=\hat PD\hat Q^{-1}$. This means that $D'=\pi(H')$, where $H'=\beta(H)$ with $\beta$ applied entry-wise. So, either $H$ is not equivalent with $H'$, but $D$ is isomorphic to $D'$ or $H$ is equivalent to $H'$, but $D$ has a larger automorphism group than $H$.

Note that there is a notion of a \emph{Hadamard equivalence} of generalized Hadamard matrices that also incorporate automorphisms of the associated group \cite{Ost22}. But we do not know if all isomorphisms of design should be of that form either. 

Doing one further step, we obtain a graph. Consider the adjacency matrix
\begin{align}\label{eq:A-original}A=\begin{pmatrix}
    0&\pi(H)\\
    \pi(H)^{\dag}&0
\end{pmatrix}.
\end{align}
Let us denote the graph associated to $A$ by $\Gamma_H$.  
Note that all the above can also be seen as a particular example of a construction of covers of the complete graph or the complete bipartite graph \cite{GH92,God96}.

Observe that any automorphism of the (order $s$) Hadamard matrix $H$, given by a pair $P,Q\in H_n^s$, induces an automorphism of $\Gamma_H$ given by the block-diagonal permutation matrix $\begin{pmatrix}\pi(P)&0\\0&\pi(Q)\end{pmatrix}$. The fact that not every automorphism of $\Gamma_H$ arises from an automorphism of $H$ as above follows from the same reasoning as  for the designs. (Besides that, the associated bipartite graph also has an extra automorphism given by exchanging the two parts. The latter operation at the level of $H$ corresponds to replacing $H$ with $H^\dag$, which does not always correspond to an automorphism of $H$.)  In the real case, one can artificially colour the two parts of $\Gamma_H$ in order to exclude this additional symmetry.  After this modification, it turns out that $\Aut(\Gamma_H) = \Aut(H)$ \cite{McK79}.  The same equality holds also for the quantum automorphism groups, as shown by the second author \cite[Proposition 6.6]{GHad}.
From the above discussion, we know that $\Aut^+(\Gamma_H) \ne \Aut^+_{H_n^{s+}}(H)$ for general $s > 2$.  In the following we show that just as in the case of classical symmetries, we in general have an embedding of quantum spaces $ \Aut^+_{H_n^{s+}}(H) \subset \Aut^+(\Gamma_H)$, in analogy with the classical situation.  In  other words:

\begin{theorem} \label{th:embedding}
Let $H$ be a  Butson matrix of size $n$ and order $s$. Then $\Aut^+_{H_n^{s+}}(H))\subseteq\Aut^+(\Gamma_H))$.
\end{theorem}  
\begin{proof}
Fix an $n \times n$ Hadamard matrix $H$ of order $s$ and let  $p$ be the fundamental representation of $\Aut^{+}_{H_n^{s+}}(H)$.  Denote $q=H^{-1}pH$.  Observe that $p,q\in H_n^{s+}$.  Denote by
$$
u_1=
\begin{pmatrix}
    a^{(0)}&a^{(1)}&\cdots&a^{(s-1)}\\
    a^{(s-1)}&a^{(0)}&\cdots&a^{(s-2)}\\
    \vdots&\vdots&\ddots&\vdots\\
    a^{(1)}&a^{(2)}&\cdots&a^{(0)}
\end{pmatrix}
\quad\text{and}\quad
u_2=
\begin{pmatrix}
    b^{(0)}&b^{(1)}&\cdots&b^{(s-1)}\\
    b^{(s-1)}&b^{(0)}&\cdots&b^{(s-2)}\\
    \vdots&\vdots&\ddots&\vdots\\
    b^{(1)}&b^{(2)}&\cdots&b^{(0)}
\end{pmatrix}
$$
the corresponding sudoku magic unitaries from Lemma~\ref{L.sudoku}.  We claim that the magic unitary
$$\begin{pmatrix}
    u_1&0\\0&u_2
\end{pmatrix}$$
commutes with the adjacency matrix $A$ of $\Gamma_H$. This will then mean that there is a surjective $*$-homomorphism $\mathcal{O}(\Aut^+(\Gamma_H))\to\mathcal{O}(\Aut^+_{H_n^{s+}}(H))$ and so that $\Aut^+_{H_n^{s+}}(H)$ equipped with the above faithful representation is a compact matrix quantum subgroup of $\Aut^+(\Gamma_H)$.

It might actually be tempting to denote the blocks above by $\pi(p)$ and $\pi(q)$, but there is a mismatch of the index ordering, which we need to resolve here anyway: The magic unitaries $u_1$ and $u_2$ are elements of $M_s(M_n(\mathcal A))$, for some (C$^\ast$-)algebra $\mathcal A$.  On the other hand the adjacency matrix $A$  of $\Gamma_H$ given by \eqref{eq:A-original} is a $2 \times 2$ anti-diagonal block matrix with non-zero blocks $\pi(H), \pi(H)^\dag  \in M_n(M_s(\mathbb C))$.  Thus to do our computations, we first must perform an index shuffle on $A$ so that its blocks are of the form $M_s(M_n(\mathbb C))$.  To do this shuffle, observe that the bipartite graph  $\Gamma_H$ naturally has vertices $\{(r_i, \omega^k)\}_{i \in [n], \ k \in \mathbb Z_s} \sqcup \{(c_j, \omega^l)\}_{j \in [n], \ l \in \mathbb Z_s}$ with the adjacency condition given by $(r_i, \omega^k) \sim (c_j, \omega^l) \iff \omega^k H_{ij} = \omega^{l} \iff H_{ij} = \omega^{l - k}$. With this in mind, we define $H^{(i)}$ to be the indicator matrix on $\omega^i$; that is, $H^{(i)}_{jk} = 1 \iff H_{jk} = \omega^i$. We may write this explicitly as $H^{(i)} = (1 + (\omega^{-i} H) + (\omega^{-i} H)^{\bullet 2} + .. + (\omega^{-i} H)^{\bullet s})$, where ${T}^{\bullet k}$ refers to the $k$th Schur power of a matrix $T$. In terms of these $H^{(k)}$, we obtain our reshuffling of $A$ as 

\begin{align} \label{eq:A-new} A = 
\begin{pmatrix}
    0&0&\cdots&0&H^{(0)}&H^{(1)}&\cdots&H^{(s-1)}\\
    0&0&\cdots&0&H^{(s-1)}&H^{(0)}&\cdots&H^{(s-2)}\\
    \vdots&\vdots&\ddots&\vdots&\vdots&\vdots&\ddots&\vdots\\
    0&0&\cdots&0&H^{(1)T}&H^{(2)T}&\cdots&H^{(0)T}\\
    H^{(0)T}&H^{(1)T}&\cdots&H^{(s-1)T}&0&0&\cdots&0\\
    H^{(s-1)T}&H^{(0)T}&\cdots&H^{(s-2)T}&0&0&\cdots&0\\
    \vdots&\vdots&\ddots&\vdots&\vdots&\vdots&\ddots&\vdots\\
    H^{(1)T}&H^{(2)T}&\cdots&H^{(0)T}&0&0&\cdots&0
\end{pmatrix}
\end{align}

We now comupute.  Succinctly writing \[u = \begin{pmatrix} u_1 & 0 \\ 0  & u_2 \end{pmatrix} \quad \& \quad A = \begin{pmatrix} 0 & X \\ X^\dag  & 0 \end{pmatrix}, \]  we obtain 
$$Au = \begin{pmatrix} 0 & u_1X \\ u_2X^\dag  & 0 \end{pmatrix} \quad \&\quad uA = \begin{pmatrix} 0 & Xu_2 \\ X^\dag u_1  & 0 \end{pmatrix}.$$  Next we observe that it suffices to to show equality of the upper right hand corners $u_1X = Xu_2$, as  the lower left hand corner relation is just the dagger of this relation. To show that $u_1X = Xu_2$, observe that since $X, X^\dag, u_2, u_2$ are (block-)circulant matrices, and products of circulant matrices are circulant,  it suffices to show that their first rows are equal.  For $0 \le i \le s-1$, we have
%
%\begin{align*}\sum_{k = 0}^{s - 1}a^{(k)H_^{(-k + i)} = \sum_{k = 0}^{s - 1}H_{\omega^{k}}q^{\omega^{-k + i}} \\
%\iff \sum_{k = 0}^{s - 1}p^{\omega^k}H_{\omega^{-k + i}} = \sum_{k = 0}^{s - 1}H_{\omega^{k - i}}q^{\omega^{-k}} \\ 
%\iff\underbrace{ \sum_{k = 0}^{s - 1}p^{\omega^k}H_{\omega^{-k + i}}}_{A(i)} = \underbrace{\sum_{k = 0}^{s - 1}H_{\omega^{-k + i}}q^{\omega^k}}_{B(i)}
%\end{align*}
%
%Then, expanding the definition of $p^{\omega^k}$ on the left hand side above, we get 

\begin{align*}
[u_1X]_{0i} &= \sum_{k = 0}^{s - 1}a^{(k)}H^{(i-k)}
             =\frac{1}{s}\sum_{k = 0}^{s - 1}\sum_{\ell = 0}^{s - 1} \omega^{-k\ell}p^{\bullet\ell}H^{(i-k)}\\
&= \frac{1}{s}\sum_{\ell = 1}^{s}p^{\bullet\ell}\sum_{k = 0}^{s - 1}\omega^{k\ell - i\ell}H^{(k)} 
= \frac{1}{s} \sum_{\ell = 1}^{s}\omega^{-i\ell}p^{\bullet \ell} H^{\bullet \ell}
\end{align*}

The last equality is justified since each $H_{\omega^{k}}$ is the indicator matrix on $\omega^{k}$ in $H$. Thus the  matrix $\sum_{k = 0}^{s - 1}\omega^{\ell k }H_{\omega^{k}}$ has entry $\omega^{\ell k}$ where $H$ has entry $\omega^k$; this is precisely the $\ell$th Schur power of $H$. A similar computation on the other side yields 

\begin{align*}
[Xu_2]_{0i} = \sum_{k = 0}^{s - 1}H^{(k)}b^{(i-k)} = \frac{1}{s} \sum_{\ell = 1}^{s}\omega^{-i\ell}H^{\bullet\ell}q^{\bullet \ell}
\end{align*}

Therefore, our desired relations are \begin{align} \label{eq:final} \frac{1}{s} \sum_{\ell = 1}^{s}\omega^{-i\ell}p^{\bullet \ell} H^{\bullet \ell} = \frac{1}{s} \sum_{\ell = 1}^{s}\omega^{-i\ell}H^{\bullet\ell}q^{\bullet \ell} \qquad (0 \le i \le s-1).
\end{align}
To prove \eqref{eq:final}, fix $1 \le \ell \le s$, and observe that 
$p^{\bullet \ell}H^{\bullet\ell} = (pH)^{\bullet\ell}$.  Indeed, from the relations in \cite[Section 2]{Bic04}, we have  for $1 \le i,j \le n$
\begin{align*}
\big((pH)^{\bullet \ell}\big)_{ij} &= \Big(\sum_{k=1}^n p_{ik}H_{kj}\Big)^\ell \\
&= \sum_{k_1, \ldots, k_\ell = 1}^n p_{ik_1}p_{ik_2} \ldots p_{ik_\ell} H_{k_1j}H_{k_2j}\ldots H_{k_\ell j} \\
&=\sum_{k_1=1}^n p_{ik_1}^\ell H_{k_1j}^\ell
\\&= (p^{\bullet \ell} H^{\bullet \ell})_{ij}.
\end{align*}
Similarly, $(Hq)^{\bullet \ell} = H^{\bullet \ell}q^{\bullet\ell}$.  Since $pH = Hq$, equality \eqref{eq:final} follows and $u$ is a quantum automorphism of $\Gamma_H$.  
%
%By universal properties of the relevant Hopf $*$-algebras, we obtain a unique unital $*$-homomorphism \[\rho:\mathcal O(\Aut^+(\Gamma_H)) \to \mathcal O(\Aut^+_{H_n^{s+}}(H)); \qquad \rho^{(2ns)}(U) = u,\] where $U$ is the fundamental representation of $\Aut^+(\Gamma_H)$ and $u$ is the magic unitary constructed above for the fundamental representation $p = [p_{ij}]$ of $\Aut^+_{H_n^{s+}}(H)$.  To see that $\rho$ is surjective, observe that 
%\[p_{ij} = \sum_{k=0}^{s-1} \omega^kp_{ij}^{\omega^k} \in \text{range}(\rho),\] and thus $\rho$ is surjective.
\end{proof}

\begin{remark} \label{rem:desymmetrized}
The proof of Theorem \ref{th:embedding} can easily be generalized {\it mutatis mutandis} to yield, for any pair of Butson matrices $H,H'$ of the same size $n$ and order $s$, a surjective $*$-homomorphism
\[\rho: \mathcal A(\text{Iso}^+(\Gamma_H, \Gamma_{H'})) \to \mathcal A_{H_n^{s+}}(H,H').\]  Here, the algebra on the left is the game algebra of the graph isomorphism game associated to $\Gamma_H,\Gamma_{H'}$.  That is the universal $*$-algebra generated by the coefficients of a magic unitary $U$ satisfying $AU = UA'$.  In particular, this gives
\begin{corollary} \label{cor:strategy-transfer}
If the Butson matrices $H,H'$ are $s$-quantum equivalent, then their associated graphs $\Gamma_H, \Gamma_{H'}$ are quantum isomorphic (in the commuting operator framework).     
\end{corollary}
\begin{proof}
Thanks to the existence of $\rho$, $\mathcal A(\text{Iso}^+(\Gamma_H, \Gamma_{H'})) \ne \{0\}$ provided $\mathcal A_{H_n^{s+}}(H,H') \ne \{0\}$, and this is sufficient to conclude thanks to \cite{BCE+}.
\end{proof}
\end{remark}

\section{Quantum symmetries of (quantum) groups} \label{sec:qsymQG}

In \cite[Question 7.1]{GHad}, the second author asked whether every fibre functor acting on the abstract diagrammatic category $\langle\Diagram{\Dmor{bcirc}2/2 (0,0.5)},\Diagram{\Dmor{circ}2/2 (0,0.5)}\rangle$ must be interpreted via (quantum) Hadamard matrices. While we keep this question open in the real case, in the complex case the answer is definitely no, because these spiders can be interpreted in terms of multiplication and comultiplication in a finite-dimensional Hopf C*-algebra as we are going to describe below. In the case when this Hopf C*-algebra comes from an abelian group, then we obtain what we had before -- quantum symmetries of the associated Fourier Hadamard matrix. But already for non-abelian groups, we obtain something new.

Another motivation for studying this particular example is the following: We observed in \cite{GHad} that the diagrammatic category $\langle\Diagram{\Dmor{bcirc}2/2 (0,0.5)},\Diagram{\Dmor{circ}2/2 (0,0.5)}\rangle$ is pure, but its natural subcategories
$$\langle\spider{2/1},\wspider{2/1}\rangle
\subseteq\langle\spider{2/1},\wspider{2/2}\rangle
\subseteq\langle\spider{2/2},\wspider{2/2}\rangle$$
are not pure. Restricting on some specific fibre functors, we may obtain extra relations that may change that.

\subsection{Hopf diagrams}
Consider a Hopf C*-algebra $\A$ with $\dim\A=n<\infty$. We denote by $\spider{2/1}$ the multiplication, by $\spider{0/1}$ the unit, by $\sqrt{n}\wspider{1/2}$ the comultiplication, and by $\frac{1}{\sqrt{n}}\wspider{1/0}$ the counit. The unique Haar state normalized such that it makes $\A$ a special Frobenius algebra is denoted by $\spider{1/0}=\spider{0/1}^\dag$. This also makes $\A$ a Hilbert space and allows one to introduce the adjoint $\dag$, which will flip every diagram about the  horizontal axis.

The axioms of Hopf algebras then mean that the white and black points satisfy all the rules of spiders. In addition, writing the axioms for Hopf algebras in terms of diagrams, we obtain the following additional rules:
\begin{equation}
\label{eq.Hopf}
\Diagram{\draw (1,0) -- (0,1);\draw (0,0) -- (1,1);
         \Dmor{circ}1/1  (0,0)
         \Dmor{circ}1/1  (1,0)
         \Dmor{bcirc}1/1 (0,1)
         \Dmor{bcirc}1/1 (1,1)}
=\frac{1}{\sqrt{n}}
\Diagram{\Dmor{bcirc}2/1 (0,0)
         \Dmor{circ}1/2  (0,1)},\quad
\Diagram{\Dmor{circ}1/2  (0,0.5) \Dmor{bcirc} 0/0 (0,0)}=\frac{1}{\sqrt n}
\Diagram{\Dmor{bcirc}0/1 (0,0.5) \Dmor{bcirc}0/1 (1,0.5)},\quad
\Diagram{\Dmor{bcirc}2/1  (0,0.5) \Dmor{circ} 0/0 (0,1)}=\frac{1}{\sqrt{n}}
\Diagram{\Dmor{circ}1/0 (0,0.5) \Dmor{circ}1/0 (1,0.5)},
\end{equation}
\begin{equation}
\label{eq.antipode}
\Diagram{\Dmor{bcirc}2/1 (0.5,2)
         \DMor{square}1/1 (0,0.5) {$\scriptstyle S$}
         \draw (1,-0.5) -- (1,1.5);
         \Dmor{circ}1/2  (0.5,-1)}
=\frac{1}{n}\Diagram{\Dmor{bcirc}0/1 (0,1)
         \Dmor{circ}1/0 (0,0)} =
\Diagram{\Dmor{bcirc}2/1 (0.5,2)
         \DMor{square}1/1 (1,0.5) {$\scriptstyle S$}
         \draw (0,-0.5) -- (0,1.5);
         \Dmor{circ}1/2  (0.5,-1)},
\end{equation}
where $S$ is the antipode.

Now, we have two bilinear forms $\spider{2/0}=\Diagram{\Dmor{bcirc}2/1 (0.5,0.5) \Dmor{bcirc}0/0 (0.5,1)}$ and $\wspider{2/0}=\Diagram{\Dmor{circ}2/1 (0.5,0.5) \Dmor{circ}0/0 (0.5,1)}$, which do not coincide. In fact, we have
\begin{equation}
\label{eq.HopfBilin}
\Diagram{\draw (-0.5,0) -- (-0.5,0.5);\draw (1.5,0.5) -- (1.5,1);
           \Dmor{bcirc}2/0 (0,1)
           \Dmor{circ}0/2  (1,0)}
=
\Diagram{\draw (-0.5,0) -- (-0.5,0.5);\draw (1.5,0.5) -- (1.5,1);
           \Dmor{circ}2/0 (0,1)
           \Dmor{bcirc}0/2  (1,0)}
=
\Diagram{\draw (1.5,0) -- (1.5,0.5);\draw (-0.5,0.5) -- (-0.5,1);
           \Dmor{circ}0/2 (0,0)
           \Dmor{bcirc}2/0  (1,1)}
=
\Diagram{\draw (1.5,0) -- (1.5,0.5);\draw (-0.5,0.5) -- (-0.5,1);
           \Dmor{bcirc}0/2 (0,0)
           \Dmor{circ}2/0  (1,1)}
=S,
\end{equation}
where the antipode in fact satisfies $S=S^*=S^\dag=S\Ctrans$.

We can introduce the dual Hilbert space $\A^*$ and arrows to the diagrams. So far, all arrows were pointing from bottom to top. Now, we can stretch out the diagrams in Eq.~\eqref{eq.HopfBilin} to express the antipode in the following way
\begin{equation}
\Diagram{\Dmor{circ}[</-] (0,0)
         \Dmor{bcirc}[>/>]  (0,1)}
=
\Diagram{\Dmor{bcirc}[</-] (0,0)
         \Dmor{circ}[>/>]  (0,1)}
=S.
\end{equation}
Using this, Eq.~\eqref{eq.antipode} gets simplified to
\begin{equation}
\label{eq.complementary}
\Diagram{\Dmor{bcirc}[->/>] (0,1)
         \Dmor{circ}[</>-] (0,0)}
=\frac{1}{n}
\Diagram{\Dmor{bcirc}[/>] (0,1)
         \Dmor{circ}[</] (0,0)} =
\Diagram{\Dmor{bcirc}[>-/>] (0,1)
         \Dmor{circ}[</->] (0,0)}
\end{equation}
So, the black and white spiders are actually \emph{complementary}.

Diagrams $\spider[>/>]$, $\spider[</<]$, $\wspider[>/>]$, and $\wspider[</<]$ provide canonical isomorphisms between $\A$ and $\A^*$. However, one has to pay attention to the fact that the black and white do not coincide and if a black isomorphism of this sort meets a white spider or vice versa, it does not merely change the direction of the arrow, but also adds the antipode.

In addition, recall that the antipode is an algebra and coalgebra antihomomorphism. This means that
$$
\Diagram{
	\Dmor{circ}[-/<]   (0,1)
	\Dmor{bcirc}[<</>] (0,0)
}
=
\Diagram{
	\Dmor{circ}[-/<]   (0,2)
    \Dmor{bcirc}[>/>]  (0,1)
	\Dmor{bcirc}[>>/-] (0,0)
    \Dmor{bcirc}[</-]  (-0.5,-1)
    \Dmor{bcirc}[</-]  (0.5,-1)
}
=
\Diagram{
	\DMor{square}[>/<]   (0,1.5) {$S$}
	\Dmor{bcirc}[>>/-] (0,0)
    \Dmor{bcirc}[</-]  (-0.5,-1)
    \Dmor{bcirc}[</-]  (0.5,-1)
}
=
\Diagram{
	\draw[->] (0,1.5) .. controls +(0,-.5) and +(0,.5) .. (1,0.5);
	\draw[->] (1,1.5) .. controls +(0,-.5) and +(0,.5) .. (0,0.5);
	\Dmor{bcirc}[>>/<] (0.5,2)
	\DMor{square}[>/-]   (0,-0.5) {$S$}
	\DMor{square}[>/-]   (1,-0.5) {$S$}
    \Dmor{bcirc}[</-]  (0,-2)
    \Dmor{bcirc}[</-]  (1,-2)
}
=
\Diagram{
	\draw[->] (0,1.5) .. controls +(0,-.5) and +(0,.5) .. (1,0.5);
	\draw[->] (1,1.5) .. controls +(0,-.5) and +(0,.5) .. (0,0.5);
	\Dmor{bcirc}[>>/<] (0.5,2)
	\Dmor{circ}[-/-]   (0,0)
	\Dmor{circ}[-/-]   (1,0)
	\Dmor{bcirc}[>/>]  (0,-1)
	\Dmor{bcirc}[>/>]  (1,-1)
    \Dmor{bcirc}[</-]  (0,-2)
    \Dmor{bcirc}[</-]  (1,-2)
}
=
\Diagram{
	\draw[->] (0,1) .. controls +(0,-.5) and +(0,.5) .. (1,0);
	\draw[->] (1,1) .. controls +(0,-.5) and +(0,.5) .. (0,0);
	\Dmor{bcirc}[>>/<] (0.5,1.5)
	\Dmor{circ}[</-]   (0,-.5)
	\Dmor{circ}[</-]   (1,-.5)
}
$$
And similarly
$$
\Diagram{
	\Dmor{circ}[-/>>]  (0,1)
	\Dmor{bcirc}[>/>]  (0,0)
}
=
\Diagram{
	\draw[->] (0,1) .. controls +(0,-.5) and +(0,.5) .. (1,0);
	\draw[->] (1,1) .. controls +(0,-.5) and +(0,.5) .. (0,0);
	\Dmor{circ}[>/--] (0.5,-.5)
	\Dmor{bcirc}[>/>] (0,1.5)
	\Dmor{bcirc}[>/>] (1,1.5)
}
$$

We will denote by $\Hopf_n$ the abstract diagrammatic category generated by the above diagrammatic relations.  Given some some diagrams $T_1,\dots,T_n\in\Hopf_n$, we will denote by $\langle T_1,\dots,T_n\rangle_n^{\Hopf}$ the subcategory of $\Hopf_n$ generated by these diagrams and by $\langle T_1,\dots,T_n\rangle_\A$ the image by the fibre functor induced by $\spider[<</>]\mapsto m_\mathcal A$, $\wspider[</>>]\mapsto\Delta_\mathcal A$.  $\Diagram{
	\DMor{square}[</>]   (0,0.5) {$S$}
	} \mapsto S_\mathcal A$.

    Below we present a few interesting applications of our diagrammatic approach to Hopf $\ast$-algebas.

\subsection{Finite quantum groups only have classical symmetries}
The first application is a short graphical proof of result that was originally  formulated in \cite{KSW15} -- establishing the absence of genuine quantum symmetries of finite quantum groups.

Given a finite-dimensional Hopf C$^*$-algebra $\A$, we can interpret it as a function algebra of a finite quantum group $G$, i.e. $\A=O(G)$. It is natural to define the {\it quantum automorphism group of $G$} to be the universal compact quantum group $H$ that acts faithfully on $G$ (i.e., $\mathcal O(H)$ co-acts on $\A=O(G)$) preserving both the multiplication as well as the comultiplication. By Tannaka-Krien duality, $H$ should correspond to the concrete category of Hilbert spaces $\langle\spider[<</>],\wspider[<</>]\rangle_\A$ and the question whether $H$ is a classical group amounts to deciding whether the tensor flip map belongs to the aforementioned category.

\begin{proposition} \label{prop:crossing-inside}
    It holds that $\Fcross\in \langle\spider[<</>],\wspider[<</>]\rangle^{\rm Hopf}_n$
\end{proposition}
\begin{proof}
In the real case, this was proven already in the (unpublished) seminal paper \cite[Example~3]{CD07} introducing ZX-calculus and complementary spiders (alternatively see \cite[Example~1]{CD09}). For the complex case, we just have to add some arrows.
$$\sqrt{n^3}
\Diagram{
	\draw [->] (-.5,-2) -- (-.5,-1);
	\draw [mid arrow] (1.5,1) -- (1.5,-1);
	\draw (-1,2) -- (-1,3);
	\Dmor{bcirc}[--/>] (0.5,2.5)
	\Dmor{circ}[-/>>]  (-0.5,1.5)
	\Dmor{circ}[--/>]  (1,1.5)
	\Dmor{bcirc}[-/>>] (0,0.5)
	\Dmor{circ}[--/>]  (0,-0.5)
	\Dmor{bcirc}[</>-] (1,-1.5)
}
=n^2
\Diagram{
	\draw [->] (-.5,-2) -- (-.5,-1);
	\draw [mid arrow] (1.5,1) -- (1.5,-1);
	\draw (-1,2) -- (-1,3);
	\draw [late arrow] (-0.5,-0.5) -- (0.5,0.5);
	\draw [late arrow] (0.5,-0.5) -- (-0.5,0.5);
	\Dmor{bcirc}[--/>] (0.5,2.5)
	\Dmor{circ}[-/>>]  (-0.5,1.5)
	\Dmor{circ}[--/>]  (1,1.5)
	\Dmor{circ}[-/>]   (-0.5,0.5)
	\Dmor{circ}[-/>]   (0.5,0.5)
	\Dmor{bcirc}[-/>]  (-0.5,-0.5)
	\Dmor{bcirc}[-/>]  (0.5,-0.5)
	\Dmor{bcirc}[</>-] (1,-1.5)
}
=n^2
\Diagram{
	\draw (-1,1) -- (-1,2);
	\draw [late arrow] (-0.5,-0.5) -- (1,0.5);
	\draw [late arrow] (1,-0.5) -- (-0.5,0.5);
	\Dmor{bcirc}[--/>] (0.5,1.5)
	\Dmor{circ}[-/>>]   (-0.5,0.5)
	\Dmor{circ}[0->/>]   (1,0.5)
	\Dmor{bcirc}[</>]  (-0.5,-0.5)
	\Dmor{bcirc}[</0>-]  (1,-0.5)
}
=n
\Diagram{
	\draw (-1,1) -- (-1,2);
	\draw [late arrow] (-0.5,-0.5) -- (1,0.5);
	\draw [late arrow] (1,-0.5) -- (-0.5,0.5);
	\Dmor{bcirc}[--/>] (0.5,1.5)
	\Dmor{circ}[-/>>]   (-0.5,0.5)
	\Dmor{circ}[/>]   (1,0.5)
	\Dmor{bcirc}[</>]  (-0.5,-0.5)
	\Dmor{bcirc}[</]  (1,-0.5)
}
=n
\Diagram{
	\draw (-.5,1) -- (-.5,2);
	\draw [mid arrow] (1,-1) .. controls +(0,1) .. (0,0.5);
	\draw [mid arrow] (0,-0.5) .. controls +(1,1) and +(1,-1) .. (1,1.5);
	\Dmor{bcirc}[-0/>] (1,1.5)
	\Dmor{circ}[-/>>]   (0,0.5)
	\Dmor{bcirc}[</>]  (0,-0.5)
}
=n
\Diagram{
	\draw (-1,1) -- (-1,2);
	\draw [mid arrow] (1,-1) .. controls +(0,.5) and +(0,-.5) .. (0,0.5);
	\draw [mid arrow] (0,-1) .. controls +(0,1) and +(1,-.5) .. (1,1.5);
	\Dmor{bcirc}[->0/>] (1,1.5)
	\Dmor{circ}[/>>-]   (0,0.5)
}
=
\Diagram{
	\draw [late arrow] (0,-1) -- (1,2);
	\draw [late arrow] (1,-1) -- (0,2);
}
$$
\end{proof}

\begin{corollary}\label{cor:noqs}
    For every finite quantum group $G$, it holds that $\Fcross\in \langle\spider[<</>],\wspider[<</>]\rangle_{O(G)}$, so it has only classical symmetries.
\end{corollary}

\subsection{Affine symmetries of finite quantum groups}
In this final section, we consider the  analogue of affine symmetries and isomorphisms of finite quantum groups. We first recall the classical notions.
Let $G$ be a finite group. We say that a permutation $\alpha\colon G\to G$ is an \emph{affine symmetry} of $G$ if $\alpha(gh^{-1}k)=\alpha(g)\alpha(h)^{-1}\alpha(k)$. We will denote the group of all affine symmetries by $\Aff G$. Likewise, we say that $G$ is \emph{affine equivalent} to $H$ if there is a map $\alpha\colon G\to H$ satisfying $\alpha(gh^{-1}k)=\alpha(g)\alpha(h)^{-1}\alpha(k)$.

We begin with an easy observation. 
\begin{lemma}
    Suppose $\alpha\colon G\to H$ is an affine equivalence. Then it is of the form $\alpha(g)=\phi(g)x$ for some group  isomorphism $\phi\colon G\to H$ and an element $x\in H$.
\end{lemma}
\begin{proof}
Let $\alpha\colon G\to H$ be an affine equivalence. Denote $x:=\alpha(e)\in H$ and define a bijection $\phi\colon G\to H$ by $\phi(g)=\alpha(g)x^{-1}$. We will show that $\phi$ is an isomorphism. Indeed,
\begin{align*}
\phi(gh^{-1})&=\alpha(gh^{-1}e)x^{-1}=\alpha(g)\alpha(h)^{-1}xx^{-1}\\&=\phi(g)xx^{-1}\phi(h)^{-1}=\phi(g)\phi(h)^{-1}.
\end{align*}
The converse is straightforward: If $\alpha(g)=\phi(g)x$, then $\alpha(gh^{-1}k)=\phi(gh^{-1}k)x=\phi(g)xx^{-1}\phi(h)^{-1}\phi(k)x=\alpha(g)\alpha(h)^{-1}\alpha(k)$.\qed
\end{proof}

As a consequence, the notion of an affine equivalence is classically not so interesting:

\begin{proposition}
Two groups are affine equivalent if and only if they are isomorphic.
\end{proposition}
\begin{proof}
By the lemma above, existence of an affine equivalence implies the existence of an isomorphism.
\end{proof}

On the other hand, one might be interested in studying the group $\Aff G=G\rtimes\Aut G$ for a given group $G$. In particular, note that taking $G=\Z_2^n$, we have $\Aff\Z_2^n=AGL(n,\Z_2)=ASL(n,\Z_2)$.

\begin{proposition}
    Let $G$ be a finite group. The representation category of $\Aff G$ is given by $\langle\Fcross,\spider[<</>],\wspider[<>/><]\rangle_{O(G)}$.
\end{proposition}
\begin{proof}
    Including the intertwiner $\Fcross$ is equivalent to saying that $\Aff G$ should be a classical group, not a genuinely quantum group. The intetwiner $\spider[<</>]$ means that $\Aff G$ should act on $G$ as a finite set. Together, the category $\langle\Fcross,\spider[<</>]\rangle_n$ corresponds to $S_n$, which is a well-known fact.

    Now, instead of $\wspider[<>/><]$, we could equivalently include $\wspider[<></>]$ or
    $$\Diagram{
        \Dmor{circ}[->-/>] (2,1)
        \Dmor{bcirc}[</-] (2,0)
        \draw [mid arrow] (1,-0.5) -- (1,0.5); \draw [mid arrow] (3,-0.5) -- (3,0.5);
    }=
    \Diagram{
        \Dmor{circ}[---/>] (2,1)
        \Dmor{circ}[>/>] (2,0)
        \Dmor{bcirc}[</-] (2,-1)
        \draw [mid arrow] (1,-1.5) -- (1,0.5); \draw [mid arrow] (3,-1.5) -- (3,0.5);
    }=
    \Diagram{
        \Dmor{circ}[---/>] (2,1)
        \DMor{square}[</>] (2,-0.5) {$S$}
        \draw [mid arrow] (1,-1.5) -- (1,0.5); \draw [mid arrow] (3,-1.5) -- (3,0.5);
    },$$
    which (up to normalization) exactly corresponds to the operation $(g,h,k)\mapsto gh^{-1}k$. So, including this intertwiner means that the permutations in $\Aff G\subset S_G$ satisfy $\alpha(gh^{-1}k)=\alpha(g)\alpha(h)^{-1}\alpha(k)$.
\end{proof}

\begin{lemma}
The following hold in $\Hopf_n$.
\begin{equation}
\label{eq.relations}
\Diagram{
	\Dmor{circ}[-/><>] (0,1)
	\Dmor{bcirc}[<</>] (0,0)
}
=
\Diagram{
	\draw [late arrow] (0.5,0) -- (0,1);
	\draw [late arrow] (1,1) -- (0.5,0);
	\draw [late arrow] (0.5,0) -- (2,1);
	\draw [late arrow] (1.5,0) -- (0,1);
	\draw [late arrow] (1,1) -- (1.5,0);
	\draw [late arrow] (1.5,0) -- (2,1);
	\Dmor{bcirc}[/>] (0,1)
	\Dmor{bcirc}[/<] (1,1)
	\Dmor{bcirc}[/>] (2,1)
	\Dmor{circ}[</]  (0.5,0)
	\Dmor{circ}[</]  (1.5,0)
},\qquad
\Diagram{
	\Dmor{circ}[-/><>] (0,1)
	\Dmor{bcirc}[>/>] (0,0)
}
=
\Diagram{
	\draw (0,0) -- (0,1);
	\draw [->] (-1,1) .. controls +(0,-.5) and (0,.5) .. (1,0);
	\draw [->] (1,1) .. controls +(0,-.5) and (0,.5) .. (-1,0);
	\Dmor{bcirc}[>/>] (-1,1.5)
	\Dmor{bcirc}[>/>] (0,1.5)
	\Dmor{bcirc}[>/>] (1,1.5)
	\Dmor{circ}[</->-] (0,-.5)
},\qquad
\Diagram{
	\Dmor{circ}[-/><>] (0,1)
	\Dmor{bcirc}[/>] (0,0)
}
=
\Diagram{
	\Dmor{bcirc}[/>] (0,0)
	\Dmor{bcirc}[/>] (1,0)
	\Dmor{bcirc}[/>] (2,0)
}
\end{equation}
\end{lemma}
\begin{proof}
    We obtain these easily from the homomorphism property of the comultiplication $\Delta$, antihomomorphism property of the antipode and unitality of the comultiplication.
\end{proof}

The quantized version of $\Aff(G)$ should, therefore, correspond to the category $\langle\spider[<</>],\wspider[<>/><]\rangle$.  More precisely, we define an abstract diagrammatic category $\AffCat_n$ to be the category consisting of all diagrams that can be constructed from $\spider[<</>]$ and $\wspider[<>/><]$ subject to the standard relations for complementary spiders and those from Eq.~\eqref{eq.relations}. This means that we have full functors $\AffCat_n\to\langle\spider[<</>],\wspider[<>/><]\rangle_n^{\rm Hopf}\to\langle\spider[<</>],\wspider[<>/><]\rangle_{O(G)}$ for any finite quantum group $G$.  We say that two quantum groups $G_1$ and $G_2$ are {\it quantum affine equivalent} if the concrete categories $\langle\spider[<</>],\wspider[<>/><]\rangle_{O(G_1)}$ and $\langle\spider[<</>],\wspider[<>/><]\rangle_{O(G_2)}$ are monoidally equivalent via a functor which sends the canonical generators to canonical generators.

We will now prove the following rather surprising result:

\begin{theorem}
\label{T.affine}
    All finite quantum groups of a fixed size are mutually quantum affine equivalent.
\end{theorem}

\begin{proof}
We will prove the result by showing that all closed diagrams in $\AffCat_n$ can be reduced to a number. This means that all fibre functors $\Aff_n \to \langle\spider[<</>],\wspider[<>/><]\rangle_n^{\rm Hopf}\to\langle\spider[<</>],\wspider[<>/><]\rangle_{O(G)}$
have the same kernel, for any $G$ of a fixed size $n$.

First, observe that every (possibly not reduced) diagram in $\AffCat_n$ satisfies that the strings going out of each white spider always have alternating directions. This is similar to the category associated to Hadamard matrices except that now it only holds for white spides, not for the black ones. The proof is the same: The property holds for the generators and is preserved under all category operations and reduction rules.

Now take any closed diagram in $\AffCat_n$. Assume it is connected. We can use rules \eqref{eq.relations} to swap white and black points. Thanks to the fact that spiders of the same colour can be merged together, this means that every closed diagram can be reduced to a single black and a white point connected by some number of edges. Since the edges have alternating directions, we can use the complementarity relation~\eqref{eq.complementary} to reduce it completely.
\end{proof}

%\begin{lemma}
%\label{L.alternating2}
%For every (possible not reduced) diagram in $\langle\spider[<</>],\spider[/>],\wspider[><></]\rangle_n^{\rm Hopf}$, the strings going out of each white spider always have alternating directions. (In particular, there always is an even number of them.)
%\end{lemma}
%\begin{proof}
%This holds for the generators. It is quite easy to see that the property is preserved under the category operations and applying the reduction rules.
%\end{proof}

%\begin{proposition}
%The category is pure. That is, all closed diagrams can be reduced to a number.
%\end{proposition}
%\begin{proof}
%Take any closed diagram in the category. Assume it is connected. We can use rules \eqref{eq.relations} to swap white and black points. Thanks to the fact that spiders of the same colour can be merged together, this means that every closed diagram, can be reduced to a single black and a white point connected by some number of edges. Since the edges have alternating directions, we can use the complementarity relation~\eqref{eq.complementary} to reduce it completely.
%\end{proof}

\begin{remark}
We finish by pointing out that there is an intersection between finite quantum groups and Hadamard matrices -- finite abelian groups / Fourier matrices. We showed in Theorem~\ref{T.weak} that all Hadamard matrices of a fixed size are mutually weakly quantum equivalent. The quantum equivalence here is supposed to generalize the classical equivalence of Hadamard matrices, where you can permute rows, permute columns, multiply rows by complex units, and multiply columns by complex units. If we focus only on Fourier type Hadamrd matrices, then Theorem~\ref{T.affine} implies a stronger result. All such matrices of a fixed size are (weakly) quantum affine equivalent, where quantum affine equivalence generalizes a similar procedure, where we are allowed to multiply only rows, but not columns by complex units. (Equivalently the other way around.)

It might be interesting to note that both weak and strong quantum equivalence quantize the ordinary equivalence, while quantum affine equivalence corresponds classically to the above described procedure which is stronger than the classical equivalence. Despite that, quantum affine equivalence does not imply strong quantum equivalence as shown by Proposition~\ref{P.Fourier}
\end{remark}

\bibliographystyle{alpha} % We choose the "plain" reference style
\bibliography{refs} % Entries are in the refs.bib file

\end{document}